\RequirePackage{ifthen}
\newboolean{MOR}
\setboolean{MOR}{false}

\ifthenelse{\boolean{MOR}}
{
\documentclass[moor]{informs4}
} 
{
\documentclass[11pt]{article}
\usepackage[margin=1in]{geometry}
}

\usepackage{graphicx}
\usepackage{color}

\newif\ifrevcolors
\revcolorsfalse 
\definecolor{darkgreen}{rgb}{0,0.45,0}
\ifrevcolors
\newcommand{\blue}[1]{\textcolor{blue}{#1}}
\newcommand{\green}[1]{\textcolor{darkgreen}{#1}}
\newcommand{\bluecolor}{\color{blue}}
\newcommand{\greencolor}{\color{darkgreen}}
\else
\newcommand{\blue}[1]{#1}
\newcommand{\green}[1]{#1}
\newcommand{\bluecolor}{}
\newcommand{\greencolor}{}
\fi
\newenvironment{greenpart}{\par\greencolor}{\par}
\newenvironment{bluepart}{\par\bluecolor}{\par}
\newcommand{\phead}[1]{\medskip\noindent\textbf{#1}\quad}

\newcommand{\abs}[1]{\lvert#1\rvert}

\usepackage{amsmath}
\usepackage{psfrag}
\usepackage{amssymb}
\usepackage{verbatim,booktabs}
\usepackage{subfigure,epsfig,url}
\usepackage{float}
\usepackage{mathrsfs}

\usepackage[makeroom]{cancel} 

\ifthenelse{\boolean{MOR}}
{
\usepackage{amsfonts}

\usepackage{natbib}
 \NatBibNumeric
 \def\newblock{\ }%
 \bibpunct[, ]{[}{]}{,}{n}{}{,}%

\usepackage{rotating}
\usepackage{fancyvrb}

\TheoremsNumberedThrough     
\ECRepeatTheorems
\JOURNAL{Mathematics of Operations Research}

\EquationsNumberedThrough    

\MANUSCRIPTNO{MOR-2026-}

\newcommand{\qedplace}{\unskip\nobreak\hfil\penalty50\hskip1em\null\nobreak\hfil
\ensuremath{\Halmos}\parfillskip=0pt\finalhyphendemerits=0\relax}
\newenvironment{prf}[1][]
{\proof{Proof.}}
{\qedplace \endproof}

\newenvironment{prfc}[1][]
{\proof{Proof #1.}}
{\qedplace \endproof}

\newtheorem{observation}{Observation}

}
{
\usepackage{hyperref}
\usepackage{amsthm}
\newcommand{\TABLE}[3]{\centering\caption{#1}#2\par\smallskip\small #3}

\newtheorem{theorem}{Theorem}
\newtheorem{lemma}{Lemma}
\newtheorem{proposition}{Proposition}
\newtheorem{corollary}{Corollary}

\newtheorem{claim}{Claim}

\newtheorem{observation}{Observation}

\newenvironment{prf}[1][]
{\begin{proof}}
{\end{proof}}

\newenvironment{prfc}[1][]
{\begin{proof}[Proof #1]}
{\end{proof}}
}

\newenvironment{cpf}
{\begin{trivlist} \item[] {\em Proof of claim. }}
{$\hfill\diamond$ \end{trivlist}}

\usepackage{xspace}

\newcommand{\floor}[1]{\lfloor#1\rfloor}
\newcommand{\ceil}[1]{\lceil#1\rceil}
\newcommand{\pare}[1]{\left(#1\right)}
\newcommand{\bra}[1]{\left\{#1\right\}}

\allowdisplaybreaks

\usepackage{enumitem} 

\DeclareMathOperator{\conv}{conv}

\newcommand{\R}{\mathbb R}
\newcommand{\Z}{\mathbb Z}
\newcommand{\Q}{\mathbb Q}
\newcommand{\N}{\mathbb N}
\newcommand{\B}{\mathcal B}
\newcommand{\norm}[1]{\lVert#1\rVert}

\newcommand{\ray}{\mathcal R}

\DeclareMathOperator{\cone}{cone}
\DeclareMathOperator{\rec}{rec{.}cone}

\renewcommand{\P}{\mathcal P}
\newcommand{\A}{\mathcal A}
\newcommand{\transp}{\mathsf T}
\newcommand{\z}{{\bf{0}}}
\renewcommand{\H}{\mathcal H}
\newcommand{\F}{\mathcal F}
\newcommand{\Qpoly}{\mathcal Q}

\newcommand{\ttl}{\blue{A} Geometric Characterization of Unbounded Integer Cubic Optimization Problems via Thin Rays}

\newcommand{\bstrct}{We study geometric characterizations of unbounded integer polynomial optimization problems.
\blue{Unboundedness along a ray characterizes unbounded integer linear and quadratic optimization problems with rational coefficients.
We show that this is no longer true in degree three, already in dimension three, in contrast with the continuous setting, where rays certify unboundedness up to degree three.
To recover a ray-based certificate, we introduce \emph{thin rays,} which are rays with an arbitrarily small neighborhood.
Our main result is that thin rays characterize unboundedness for integer cubic optimization problems over arbitrary rational polyhedra, in every dimension.
As a special case, we obtain a characterization of unbounded integer quadratic optimization problems with irrational coefficients, a setting in which thin rays are again necessary.
Degree three is also the threshold for certificates of this type: from degree four on, no ray certifies unboundedness, even if we allow a neighborhood of any width around it, already in dimension two and with rational coefficients.
Together, these results give a complete picture of when unboundedness of an integer polynomial optimization problem is certified by a ray, by a thin ray, or by neither.}}

\newcommand{\kywrds}{integer polynomial optimization;
integer cubic optimization;
unbounded problems;
unbounded rays}

\newcommand{\MSCcodes}{Primary: 90C10; secondary: 90C26}

\ifthenelse{\boolean{MOR}}
{
\TITLE{\ttl}
\RUNTITLE{\ttl}

\ABSTRACT{\bstrct}
}
{
\title{\ttl}

\author{Alberto Del Pia
\thanks{Department of Industrial and Systems Engineering \& Wisconsin Institute for Discovery,
University of Wisconsin-Madison.
E-mail: {\tt delpia@wisc.edu}}
}
}

\date{September 21, 2026}

\begin{document}

\maketitle

\ifthenelse{\boolean{MOR}}
{}
{
\begin{abstract}
\bstrct
\end{abstract}



\bigskip

\noindent
\hangindent=2cm
\emph{Key words:}
\kywrds

\bigskip

\noindent
\emph{Mathematics Subject Classification:} 
\MSCcodes
}

\section{Introduction}
\label{sec intro}

A central question in optimization is to characterize the conditions under which a problem becomes unbounded.
For both integer linear and integer quadratic optimization \blue{with rational data}, this question admits an elegant geometric answer:
informally, a problem is unbounded if and only if it is unbounded along some ray.
This result has deep geometric significance: it shows that whenever the objective diverges on the feasible points, it does so in the simplest possible manner, along a single direction.

We now introduce some notation to state this result precisely.
A \emph{polyhedron} is a set of the form $\P = \bra{x \in \R^n : Ax \le b}$, where $A \in \R^{m \times n}$ and $b \in \R^m$.
We say that $\P$ is rational if $A \in \Q^{m \times n}$ and $b \in \Q^m$.
A \emph{ray} in $\R^n$ is a set of the form $\ray(y, d) = \{ y + \lambda d : \lambda \ge 0\}$, for some $y \in \R^n$ and some nonzero $d \in \R^n$.
A \emph{ray} of a polyhedron $\P$ is defined as a ray $\ray(y, d)$ fully contained in $\P$, and this happens if and only if $y \in \P$, and $d$ is in the recession cone of $\P$, which we denote by $\rec(\P)$ (see, e.g., \cite{SchBookIP}).
\blue{The problems that we study in the first part of this introduction are integer polynomial optimization problems with rational data, that is, problems of the form
\begin{align}
\label{eq ipo}
\inf \bra{f(x) : x \in \P \cap \Z^n},
\end{align}
where $\P \subseteq \R^n$ is a rational polyhedron and $f \colon \R^n \to \R$ is a polynomial function with rational coefficients.
Integer linear, quadratic, and cubic optimization are the special cases of problem~\eqref{eq ipo} in which $f$ has degree at most one, two, and three, respectively.
We say that problem~\eqref{eq ipo} is \emph{unbounded} if $f$ is unbounded below on $\P \cap \Z^n$.}
We can now formally state this known characterization for integer linear and quadratic optimization.
We refer the reader to \blue{theorem~5.1 in~\cite{Mey74}} for the linear case and to theorem~4 in~\cite{dPDeyMol17MPA} for the quadratic case.

\begin{theorem}
\label{th known}
Let $\P \subseteq \R^n$ be a rational polyhedron and let $f \colon \R^n \to \R$ be a linear or quadratic function with rational coefficients.
Then $f$ is unbounded below on $\P \cap \Z^n$, if and only if there exists a ray $\ray(y, d)$ of $\P$ such that $f$ is unbounded below on $\ray(y, d) \cap \Z^n$.
\end{theorem}

\blue{In other words, Theorem~\ref{th known} states that a problem~\eqref{eq ipo} of degree at most two is unbounded if and only if it is unbounded along some ray of $\P$.
It is natural to ask whether the same holds for problems~\eqref{eq ipo} of higher degree.
It is simple to see that the answer is negative for problems~\eqref{eq ipo} of degree four or higher, and in a strong sense: already in dimension two, unboundedness may be certified by a parabola, but by no ray, even if we allow a neighborhood of any width around the ray.
In the statement, we denote by $\B_\epsilon$ the ball in $\R^n$ centered at the origin with radius $\epsilon$, i.e., $\B_\epsilon := \bra{x \in \R^n : \norm{x} \le \epsilon}$.}

\begin{proposition}
\label{prop quartic example}
\blue{Let $k \ge 4$ be an integer.
There exist a rational polyhedron $\P \subseteq \R^2$ and a polynomial function $f \colon \R^2 \to \R$ of degree $k$ with rational coefficients, such that $f$ is unbounded below on $\P \cap \Z^2$, and for every ray $\ray(y,d)$ of $\P$ and every constant $\epsilon > 0$, $f$ is bounded below on $\P \cap \pare{\ray(y,d) + \B_\epsilon}$.}
\end{proposition}

\blue{The proof of Proposition~\ref{prop quartic example} builds on a classical example of Frank and Wolfe~\cite{FraWol56}.
For problems~\eqref{eq ipo} of degree four or higher, any geometric characterization of unboundedness must therefore involve more complex objects than rays, such as polynomial curves.}

\blue{This leaves open the case of problems~\eqref{eq ipo} of degree three, which is the only one in which the answer to the above question is not settled by Theorem~\ref{th known} or by Proposition~\ref{prop quartic example}.
Here rays are known to suffice only for $n \le 2$, by theorem~1.3 in~\cite{dPHilWeiZem16MOR}, while for every $n \ge 3$ the question was open.
It is instructive to first look at the continuous setting, that is, at the continuous relaxation $\inf \bra{f(x) : x \in \P}$ of problem~\eqref{eq ipo}, in which $f$ is minimized over all points of $\P$, rather than over its integer points.
In the continuous setting, an analogue of Theorem~\ref{th known} is valid up to degree three: for degree one, this is a basic fact of linear programming (see, e.g., \cite{SchBookIP}), while for degree two and three it was proved in~\cite{Vav90} and in~\cite{BiedPHil23MPB}, respectively.
For polynomials of degree four or higher, instead, rays do not certify unboundedness even in the continuous setting: this is proposition~5.3 in~\cite{BiedPHil23MPB}, and it also follows from Proposition~\ref{prop quartic example}, where boundedness holds over the reals, and on neighborhoods of rays of any width.}

\blue{We first show that, in sharp contrast with the continuous setting, Theorem~\ref{th known} does not extend to problems~\eqref{eq ipo} of degree three, already for $n=3$.}

\begin{proposition}
\label{prop cubic example}
There exist a rational polyhedron $\P \subseteq \R^3$ and a cubic polynomial $f \colon \R^3 \to \R$ with rational coefficients, such that $f$ is unbounded below on $\P \cap \Z^3$, and for every ray $\ray(y, d)$ of $\P$, $f$ is bounded below on $\ray(y, d) \cap \Z^3$.
\end{proposition}

The failure of Theorem~\ref{th known} for problems~\eqref{eq ipo} of degree three raises a fundamental question:
\begin{quote}
\itshape Can we give a geometric characterization of unbounded integer cubic optimization problems?
\end{quote}
\blue{In the continuous setting, we have seen that rays certify unboundedness for problems of degree three.
In the integer setting, Proposition~\ref{prop cubic example} rules them out, but for a reason of a specific kind: the direction of the ray along which $f$ is unbounded below may be irrational, and a ray with irrational direction contains at most one integer point.
This suggests that rays might still suffice, provided we allow a neighborhood around each ray.
Our main result shows that this is indeed the case.}

\begin{theorem}
\label{th main}
Let $\P \subseteq \R^n$ be a rational polyhedron and let $f \colon \R^n \to \R$ be a \blue{polynomial function of degree at most three}.
Then $f$ is unbounded below on $\P \cap \Z^n$, if and only if there exists a ray $\ray(y, d)$ of $\P$ such that, for every constant $\epsilon > 0$, $f$ is unbounded below on $\P \cap \Z^n \cap \pare{\ray(y,d) + \B_\epsilon}$.
\end{theorem}

\blue{We stress that the neighborhood in Theorem~\ref{th main} can be taken arbitrarily small, since unboundedness is required for every $\epsilon > 0$, and not only for some fixed $\epsilon$.
We refer to rays with an arbitrarily small neighborhood as ``thin rays.''
Theorem~\ref{th main} is tight in two respects: by Proposition~\ref{prop cubic example} it fails with $\epsilon = 0$, and by Proposition~\ref{prop quartic example} it does not extend to problems~\eqref{eq ipo} of degree four, already for $n=2$, where not even thin rays certify unboundedness.
Together with Theorem~\ref{th known} and the continuous results recalled above, we obtain the complete picture for problem~\eqref{eq ipo} summarized in Table~\ref{tab rational}: the integer and the continuous setting agree except in degree three, where rays must be replaced by thin rays.}

\begin{table}[htbp]
\TABLE
{\blue{Certificates of unboundedness for problem~\eqref{eq ipo} and for its continuous relaxation, by degree $k$ of $f$.}\label{tab rational}}
{\bluecolor
\begin{tabular}{@{}lll@{}}
\toprule
Degree $k$ & Continuous setting & Integer setting \\
\midrule
$1$ & rays~\cite{SchBookIP} & rays~\cite{Mey74} \\
$2$ & rays~\cite{Vav90} & rays~\cite{dPDeyMol17MPA} \\
$3$ & rays~\cite{BiedPHil23MPB} & thin rays (Prop.~\ref{prop cubic example}, Thm.~\ref{th main}) \\
$\ge 4$ & none~\cite{BiedPHil23MPB} (Prop.~\ref{prop quartic example}) & none (Prop.~\ref{prop quartic example}) \\
\bottomrule
\end{tabular}}
{\blue{``Rays'': the problem is unbounded if and only if it is unbounded along some ray of $\P$. ``Thin rays'': rays do not suffice, but thin rays do. ``None'': neither rays nor thin rays suffice.}}
\end{table}

\blue{Theorem~\ref{th main} is, to the best of our knowledge, the first structural result for integer cubic optimization.}
This stands in stark contrast with the linear case, which has been extensively studied (see, e.g., \cite{SchBookIP}), and with the quadratic case, which has attracted significant attention \blue{(see~\cite{Bal95,dP18MPB,dP19SIOPT,dP23bMPA,dP25xSIOPT,dPDeyMol17MPA,dPMa22MPA,dPPos19MPA,dPWei14SODA,EibGanKnoOrd19,HilWeiZem16,LeeOnnRomWei12,Lok17} and the references therein)}\blue{, while for higher degree the available results are essentially algorithmic, such as the FPTAS in fixed dimension of~\cite{DeLHemKopWei06b,DeLHemKopWei08}}.
\blue{Thus, Theorem~\ref{th main} marks an important step toward a broader geometric understanding of integer polynomial optimization beyond the quadratic case.}

\blue{The ``if'' direction of Theorem~\ref{th main} is immediate.
We now outline the proof of the ``only if'' direction, which is an iterated projection argument.
Starting from a sequence of integer points of $\P$ along which $f$ tends to $-\infty$, we successively project the sequence, along suitable limit directions, onto faces of the recession cone of $\P$ of smaller and smaller dimension.
The projected points are in general not integer, and the main technical point is to show that they nevertheless inherit the unboundedness of $f$.
Each projection strictly decreases the dimension of the face, so the process cannot continue indefinitely; and it stops only when a ray of $\P$ starting at an integer point along which $f$ is unbounded below has been found.
The analysis rests on the interplay between the tensor, matrix, and vector coefficients of $f$ along the projection directions: unlike the linear and quadratic cases, where vectors and matrices suffice, the cubic setting requires tensors, and we use tensor notation, in particular tensor contractions, throughout the paper.}

\blue{In the remainder of this introduction, we drop the rationality assumption on the coefficients of the objective function, and we consider problems of the form
\begin{align}
\label{eq ipo real}
\inf \bra{f(x) : x \in \P \cap \Z^n},
\end{align}
where $\P \subseteq \R^n$ is still a rational polyhedron, but $f \colon \R^n \to \R$ is now a polynomial function with arbitrary real coefficients.
As before, problem~\eqref{eq ipo real} is \emph{unbounded} if $f$ is unbounded below on $\P \cap \Z^n$.}
\blue{Theorem~\ref{th known} assumes rationality, whereas Theorem~\ref{th main} does not, and therefore applies to problems~\eqref{eq ipo real} of degree at most three.
This is of interest already in degree two, where rays are no longer sufficient, even for $n=2$.}

\begin{proposition}
\label{prop quadratic example}
There exist a rational polyhedron $\P \subseteq \R^2$ and a quadratic polynomial $f \colon \R^2 \to \R$ \blue{with real coefficients}, such that $f$ is unbounded below on $\P \cap \Z^2$, and for every ray $\ray(y, d)$ of $\P$, $f$ is bounded below on $\ray(y, d) \cap \Z^2$.
\end{proposition}

\blue{Thin rays, instead, do suffice: the special case of Theorem~\ref{th main} in which $f$ is quadratic reads as follows.}

\begin{corollary}
\label{th quad}
Let $\P \subseteq \R^n$ be a rational polyhedron and let $f \colon \R^n \to \R$ be a quadratic function.
Then $f$ is unbounded below on $\P \cap \Z^n$, if and only if there exists a ray $\ray(y, d)$ of $\P$ such that, for every constant $\epsilon > 0$, $f$ is unbounded below on $\P \cap \Z^n \cap \pare{\ray(y,d) + \B_\epsilon}$.
\end{corollary}

\blue{For rational coefficients, Corollary~\ref{th quad} follows immediately from Theorem~\ref{th known}, so its content lies in the objective functions with irrational coefficients.
For degree one, instead, Theorem~\ref{th known} remains valid for real coefficients, by standard arguments in integer linear optimization.
Table~\ref{tab real} summarizes the picture: in the continuous setting the rationality of the data plays no role, and the certificates are the same as in Table~\ref{tab rational}, while in the integer setting the only difference from Table~\ref{tab rational} is in degree two, where rays must be replaced by thin rays.}

\begin{table}[htbp]
\TABLE
{\blue{Certificates of unboundedness for problem~\eqref{eq ipo real} and for its continuous relaxation, by degree $k$ of $f$.}\label{tab real}}
{\bluecolor
\begin{tabular}{@{}lll@{}}
\toprule
Degree $k$ & Continuous setting & Integer setting \\
\midrule
$1$ & rays~\cite{SchBookIP} & rays~\cite{Mey74} \\
$2$ & rays~\cite{BiedPHil23MPB} & thin rays (Prop.~\ref{prop quadratic example}, Cor.~\ref{th quad}) \\
$3$ & rays~\cite{BiedPHil23MPB} & thin rays (Prop.~\ref{prop cubic example}, Thm.~\ref{th main}) \\
$\ge 4$ & none~\cite{BiedPHil23MPB} (Prop.~\ref{prop quartic example}) & none (Prop.~\ref{prop quartic example}) \\
\bottomrule
\end{tabular}}
{\blue{Notation as in Table~\ref{tab rational}.}}
\end{table}

This paper is organized as follows. Section~\ref{sec notation} introduces the notation for linear, bilinear, and trilinear forms\blue{, together with the contractions and norms that we use}.
\blue{In Section~\ref{sec main}, we reduce the ``only if'' direction of Theorem~\ref{th main} to two propositions: Proposition~\ref{prop int ray}, which produces a ray of $\P$ starting at an integer point along which $f$ is unbounded below, and Proposition~\ref{prop thin ray}, which turns that ray into a thin ray.
That section also contains the proof of Proposition~\ref{prop thin ray}.
Section~\ref{sec other} contains the proofs of Propositions~\ref{prop quartic example}, \ref{prop cubic example}, and~\ref{prop quadratic example}, the last two of which rely on Proposition~\ref{prop thin ray}.
Section~\ref{sec proof cor quad} gives a direct proof of Corollary~\ref{th quad}, which serves as a warm-up for the proof of Proposition~\ref{prop int ray}.
Proposition~\ref{prop int ray}, which is the heart of the paper and completes the proof of Theorem~\ref{th main}, is proved in Section~\ref{proof prop int ray}, which opens with an overview of the argument.}

\section{\blue{Notation}}
\label{sec notation}

In this section, we present our notation for linear, bilinear, and trilinear forms, which will be used throughout the paper.
We denote by $\N$ the set $\Z_{\ge 1}$ of positive integers.

Given a vector $V \in \R^n$ and $x \in \R^n$, we denote by $V[x]$ the linear form
\[
V[x] = \sum_{i=1}^n V_i x_i.
\]
Given a matrix $M \in \R^{n \times n}$ and vectors $x,y \in \R^n$, we denote by $M[x,y]$ the bilinear form
\[
M[x,y] = \sum_{i,j=1}^n M_{ij} x_i y_j.
\]
Given a tensor $T \in \R^{n \times n \times n}$ and vectors $x,y,z \in \R^n$, we denote by $T[x,y,z]$ the trilinear form
\[
T[x,y,z] = \sum_{i,j,k=1}^n T_{ijk} x_i y_j z_k.
\]
If $M$ is symmetric, i.e., $M_{ij} = M_{ji}$ for all $i,j$, then $M[x,y] = M[y,x]$.  
Similarly, if $T$ is symmetric, meaning its entries are invariant under any permutation of the indices, then $T[x,y,z]$ is invariant under any permutation of its arguments.

With this notation, every polynomial function of degree at most three can be written compactly.

%
%




\begin{observation}
\label{obs f notation}
Let $f : \R^n \to \R$ be a polynomial function of degree at most three.
Then, there exist a symmetric tensor $T \in \R^{n \times n \times n}$, a symmetric matrix $M \in \R^{n \times n}$, a vector $V \in \R^n$, and a constant $c \in \R$ such that
\[
f(x) = T[x,x,x] + M[x,x] + V[x] + c.
\]
Moreover, if $f$ has rational coefficients, $T$, $M$, $V$, and $c$ can be taken rational.
\end{observation}

\begin{prf}
Let $f : \R^n \to \R$ be a polynomial of degree at most three. Then $f$ can be written as
\[
f(x) = \sum_{i,j,k=1}^n \alpha_{ijk} x_i x_j x_k 
      + \sum_{i,j=1}^n \beta_{ij} x_i x_j 
      + \sum_{i=1}^n \gamma_i x_i 
      + c,
\]
for some coefficients $\alpha_{ijk}, \beta_{ij}, \gamma_i \in \R$ and constant $c \in \R$.

We define a symmetric tensor $T \in \R^{n \times n \times n}$ by symmetrizing the coefficients $\alpha_{ijk}$:
\[
T_{ijk} := \frac{1}{6} (\alpha_{ijk} + \alpha_{ikj} + \alpha_{jik} + \alpha_{jki} + \alpha_{kij} + \alpha_{kji}) \,.
\]
Then $T$ is symmetric, and for all $x \in \R^n$, we have
\[
\sum_{i,j,k=1}^n \alpha_{ijk} x_i x_j x_k = \sum_{i,j,k=1}^n T_{ijk} x_i x_j x_k = T[x,x,x].
\]

Similarly, define a symmetric matrix $M \in \R^{n \times n}$ by
\[
M_{ij} := \frac{1}{2} (\beta_{ij} + \beta_{ji}) \,,
\]
so that
\[
\sum_{i,j=1}^n \beta_{ij} x_i x_j = \sum_{i,j=1}^n M_{ij} x_i x_j = M[x,x].
\]

Finally, set $V = (\gamma_1,\dots,\gamma_n)^\transp \in \R^n$. Then
\[
\sum_{i=1}^n \gamma_i x_i = V[x] \,.
\]

Putting everything together, we obtain
\[
f(x) = T[x,x,x] + M[x,x] + V[x] + c.
\]
\blue{If, moreover, $f$ has rational coefficients, then $\alpha_{ijk}$, $\beta_{ij}$, $\gamma_i$, and $c$ can be taken rational, and so are the entries of $T$, $M$, and $V$, since they are averages of the $\alpha_{ijk}$ and of the $\beta_{ij}$, and coincide with the $\gamma_i$, respectively.}
\end{prf}


Next, we define standard contractions in matrix and tensor calculus.
Let $M \in \R^{n \times n}$ be a matrix and let $d \in \R^n$. 
We denote by $M[d] \in \R^n$ the \emph{contraction} of $M$ with the vector $d$ along its second index:
\[
(M[d])_i := \sum_{j=1}^n M_{ij} d_j.
\]
Contraction along the first index is defined analogously.
Tensor contractions are defined similarly.
Let $T \in \R^{n \times n \times n}$ be a tensor and let \blue{$u,v \in \R^n$}. 
We denote by \blue{$T[u,v] \in \R^n$} the \emph{contraction} of $T$ with \blue{$u$ and $v$} along its second and third indices:
\[
{\bluecolor (T[u,v])_i := \sum_{j,k=1}^n T_{ijk} u_j v_k.}
\]
Contraction along any other pair of indices is defined similarly.
If $M$ is symmetric, then contractions along the first and second indices coincide.  
Similarly, if $T$ is symmetric, all contractions coincide.

We conclude this section by defining the norms that we will use.
Given a vector $x \in \R^n$, we denote by $\norm{x}$ its Euclidean norm.
\blue{Given a matrix $M \in \R^{m \times n}$, we denote by $\norm{M}$ its \emph{operator norm}
\[
\norm{M} := \sup_{\substack{x \in \R^{n} \\ \norm{x}=1}} \norm{M x},
\]
so that $\norm{Mx} \le \norm{M} \, \norm{x}$ for every $x \in \R^{n}$.}
Given \blue{a square matrix $M \in \R^{n\times n}$} and a tensor $T\in\R^{n\times n\times n}$, we define their \emph{spectral norms} by
\begin{align*}
\norm{M} & := \sup_{\substack{x,y \in \R^n \\ \norm{x}=\norm{y}=1}} |M[x,y]|, \\
\norm{T} & := \sup_{\substack{x,y,z \in \R^n \\ \norm{x}=\norm{y}=\norm{z}=1}} |T[x,y,z]|.
\end{align*}
\blue{For a square matrix the two definitions agree, so the notation is unambiguous.
By homogeneity, $\abs{M[x,y]} \le \norm{M} \, \norm{x} \, \norm{y}$ and $\abs{T[x,y,z]} \le \norm{T} \, \norm{x} \, \norm{y} \, \norm{z}$, for every $x,y,z \in \R^n$, while $\abs{V[x]} \le \norm{V} \, \norm{x}$ for every $V,x \in \R^n$, by the Cauchy--Schwarz inequality.}

\section{\blue{From rays to thin rays}}
\label{sec main}
\label{sec proof prop thin ray}

To establish Theorem~\ref{th main}, it suffices to prove the following two results:


\begin{proposition}
\label{prop int ray}
Let $\P$ be a rational polyhedron in $\R^n$, and let $f : \R^n \to \R$ be a polynomial function of degree at most three that is unbounded below on $\P \cap \Z^n$.
Then, there exists a ray $\ray(y,d)$ of $\P$ with $y \in \Z^n$ such that $f$ is unbounded below on $\ray(y,d)$.
\end{proposition}

\begin{proposition}
\label{prop thin ray}
Let $\P$ be a rational polyhedron in $\R^n$, and let $f : \R^n \to \R$ be a polynomial function of degree at most three.
Assume that there exists a ray $\ray(y,d)$ of $\P$ with $y \in \Z^n$ such that $f$ is unbounded below on $\ray(y,d)$.
Then, \blue{there exists $\bar y \in \P \cap \Z^n$ such that,} for every constant $\epsilon > 0$, $f$ is unbounded below on $\P \cap \Z^n \cap \pare{\ray(\bar y,d) + \B_\epsilon}$.
\end{proposition}

\blue{Indeed, the ``if'' direction of Theorem~\ref{th main} is trivial, since $\P \cap \Z^n \cap \pare{\ray(y,d) + \B_\epsilon} \subseteq \P \cap \Z^n$, while the ``only if'' direction follows by applying first Proposition~\ref{prop int ray} and then Proposition~\ref{prop thin ray}.
The two propositions play very different roles.
Proposition~\ref{prop int ray} is the core of the paper: it states that unboundedness on the integer points of $\P$ is always witnessed by a ray of $\P$ starting at an integer point.
Proposition~\ref{prop thin ray} then shows that such a ray can be turned into a thin ray with the same direction, possibly after changing its starting point, even when its direction is irrational.
In this section, we prove Proposition~\ref{prop thin ray}, for which we need some preparation.}



We will rely on the following standard result; see, for example, proposition~1 in~\cite{dP25xSIOPT}.

\begin{lemma}
\label{lem projection}
Let $\A$ be a rational affine subspace of $\R^n$ of dimension $n'$ with $\A \cap \Z^n \neq \emptyset$.
Then, there exists a map $\pi : \R^{n'} \to \R^n$ of the form $\pi(x') = \tilde x + \blue{U} x'$ with $\tilde x \in \Z^n$ and $\blue{U} \in \Q^{n \times n'}$ of full rank, such that 
\[
\begin{aligned}
\A & = \pi(\R^{n'}), \\
\A \cap \Z^n & = \pi(\Z^{n'}).
\end{aligned}
\]
\end{lemma}



\blue{We also rely on} the following result, which uses Dirichlet's simultaneous approximation theorem for rational polyhedra.
For a discussion of the continued fraction method, see, for example, section~6.1 in~\cite{SchBookIP}.

\begin{lemma}
\label{lem Dirichlet poly}
Let $\P$ be a rational polyhedron in $\R^n$, and let $\ray(y,d)$ be a ray of $\P$ with $y \in \Z^n$.
Then, for every $\epsilon > 0$ and every $\bar \lambda > 0$, there exists a vector in $\P \cap \Z^n$ at distance less than $\epsilon$ from the half-line $\bra{y + \lambda d : \lambda \ge \bar \lambda}$.
\end{lemma}


\begin{prf}
Fix $\epsilon > 0$ and $\bar \lambda > 0$.
Let $\P_I$ denote the integer hull of $\P$, i.e., $\P_I = \conv(\P \cap \Z^n)$.
It is well known that $\P_I$ is a rational polyhedron with $\rec(\P_I) = \rec(\P)$\blue{; see theorem~16.1 in~\cite{SchBookIP}}.
It follows that $\ray(y,d)$ is a ray of $\P_I$.

Let $\F$ be the smallest face of $\P_I$ containing $\ray(y,d)$.
\blue{If $\F$ has dimension one, then, since $\F$ is a face of the rational polyhedron $\P_I$ that contains $\ray(y,d)$, the vector $d$ is a positive multiple of an integer vector $u$ parallel to $\F$.
As $y \in \Z^n$, it follows that $y + k u \in \F \cap \Z^n \subseteq \P \cap \Z^n$, for every nonnegative integer $k$.
Hence, the half-line $\bra{y + \lambda d : \lambda \ge \bar \lambda}$ contains infinitely many vectors in $\P \cap \Z^n$, and we are done.}
So we now assume that $\F$ has dimension at least two.
Let $\A$ be the smallest affine subspace of $\R^n$ containing $\F$.
Observe that each vector in the set $\bra{y + \lambda d : \lambda > 0}$ is in the relative interior of $\F$\blue{: otherwise, one such vector would lie in a proper face of $\F$, which would then contain $\ray(y,d)$, contradicting the minimality of $\F$}.

\blue{We claim that, by reducing $\epsilon$ if necessary, we can assume that the vectors in $\A$ at distance less than $\epsilon$ from $\bra{y + \lambda d : \lambda \ge \bar \lambda}$ are contained in $\F$.
If $\F = \A$ there is nothing to prove, so we assume $\F \neq \A$, and we write $\F = \bra{x \in \A : g_i^\transp x \le h_i, \ i \in \bra{1,\dots,r}}$, where $r \ge 1$, $g_i \neq \z$ for every $i$, and no inequality is satisfied with equality by every point of $\F$, so that the relative interior of $\F$ is the set of vectors in $\A$ that satisfy all the inequalities strictly.
Since $d \in \rec(\F)$, we have $g_i^\transp d \le 0$, thus $h_i - g_i^\transp (y + \lambda d) \ge h_i - g_i^\transp (y + \bar \lambda d) > 0$ for every $\lambda \ge \bar \lambda$ and $i \in \bra{1,\dots,r}$, where the last inequality holds because $y + \bar \lambda d$ is in the relative interior of $\F$.
Hence every vector in $\A$ at distance less than $\min_{i} \frac{h_i - g_i^\transp (y + \bar \lambda d)}{\norm{g_i}}$ from the half-line satisfies all the inequalities, and the claim follows.
Reducing $\epsilon$ is without loss of generality, since the statement of the lemma for a smaller $\epsilon$ implies it for the original one.}

If $\F$ has dimension $n$, then $\F = \P_I$ and $\A = \R^n$.
It follows from Dirichlet's simultaneous approximation theorem (see, e.g., lemma~2.2 in \cite{BasConCorZam10b}) that there exists a vector $z \in \Z^n$ at distance less than $\epsilon$ from the half-line $\bra{y + \lambda d : \lambda \ge \bar \lambda}$.
We then obtain $z \in \F = \P_I \subseteq \P$.
Therefore, in the remainder of the proof, we assume that $\F$ has dimension at most $n-1$.

Denote by $n'$ the dimension of $\A$, which coincides with the dimension of $\F$, and is therefore at most $n-1$.
Since $\P_I$ is rational, there exist $W \in \Q^{l \times n}$ and $w \in \Q^{l}$ such that $\A = \bra{x \in \R^n : Wx = w}$.
Since $\A$ contains integer points, it follows from Lemma~\ref{lem projection} that there exists a map $\pi : \R^{n'} \to \R^n$ of the form $\pi(x') = \tilde x + \blue{U} x'$ with $\tilde x \in \Z^n$ and $\blue{U} \in \Q^{n \times n'}$ of full rank, such that 
\[
\begin{aligned}
\A & = \pi(\R^{n'}), \\
\A \cap \Z^n & = \pi(\Z^{n'}).
\end{aligned}
\]
Let $y' := \pi^{\leftarrow}(y)$ and $d' := \pi^{\leftarrow}(y+d)-y'$, where $\pi^{\leftarrow}$ denotes the inverse of $\pi$, so that $\pi(y')=y$ and $\pi(y'+d')=y+d$.
\blue{Note that $y' \in \Z^{n'}$, since $y \in \A \cap \Z^n = \pi(\Z^{n'})$, and that $d' \neq \z$, since $\blue{U} d' = d \neq \z$.}
Let $\epsilon' := \epsilon / \norm{\blue{U}}$\blue{, which is well defined since $U$ has full rank, hence $\norm{U} > 0$}.
It follows from Dirichlet's simultaneous approximation theorem (see, e.g., lemma~2.2 in \cite{BasConCorZam10b}) that there exists a vector $z' \in \Z^{n'}$ at distance less than $\epsilon'$ from the half-line $\bra{y' + \lambda d' : \lambda \ge \bar \lambda}$.
Let $z := \pi(z')$.
Then, $z \in \A \cap \Z^n$.

We show that $z$ is at distance less than $\epsilon$ from the half-line $\bra{y + \lambda d : \lambda \ge \bar \lambda}$.
To see this, let $q'$ be a vector in the half-line $\bra{y' + \lambda d' : \lambda \ge \bar \lambda}$ such that $\norm{z'-q'} < \epsilon'$.
Define $q := \tilde x + \blue{U} q'$, which is in the half-line $\bra{y + \lambda d : \lambda \ge \bar \lambda}$.
Then,
\[
z-q = (\tilde x + {\bluecolor U}z') - (\tilde x + {\bluecolor U}q') = {\bluecolor U}(z'-q'),
\]
thus
\[
\norm{z-q} = \norm{{\bluecolor U}(z'-q')}
\le \norm{{\bluecolor U}} \, \norm{z'-q'}
< \norm{{\bluecolor U}} \,\epsilon'
{\bluecolor = \epsilon.}
\]
\blue{Since $z \in \A$ and $z$ is at distance less than $\epsilon$ from the half-line, it follows that $z \in \F \subseteq \P_I \subseteq \P$, and therefore $z \in \P \cap \Z^n$.}
\end{prf}

%
%



%


Given a function $f : \R^n \to \R$ and a ray $\ray(y,d)$ in $\R^n$, the \emph{restriction} of $f$ to $\ray(y,d)$ is the univariate function $f_{\ray(y,d)} : \R_{\ge 0} \to \R$ defined by
\[
f_{\ray(y,d)}(\lambda) := f(y+\lambda d).
\]
\blue{If $f$ is a polynomial function, then $f_{\ray(y,d)}$ is a univariate polynomial, and $f$ is unbounded below on $\ray(y,d)$ if and only if $\lim_{\lambda \to + \infty} f_{\ray(y,d)}(\lambda) = -\infty$.}

We are now ready to present our proof of Proposition~\ref{prop thin ray}.
\blue{The idea is the following.
We compare the restriction of $f$ to $\ray(y,d)$, which is a univariate polynomial of degree at most three, with its restrictions to the nearby rays $\ray(y+z,d)$ with $\norm{z} \le \epsilon$: the leading coefficient $T[d,d,d]$ is unchanged, while the remaining coefficients are perturbed by quantities that vanish as $\epsilon$ tends to zero.
Whenever the coefficient responsible for the unboundedness of $f$ on $\ray(y,d)$ survives such a perturbation, $f$ tends to $-\infty$ on the entire neighborhood of the ray, and Lemma~\ref{lem Dirichlet poly} supplies the integer points of $\P$ that are needed to conclude.
The delicate case is the one in which the restriction of $f$ to $\ray(y,d)$ is linear: there the coefficient of $\lambda^2$ is zero and is perturbed by $3 z^\transp T[d,d]$, which may be positive, so that $f$ tends to $+\infty$ on some of the nearby rays.
If $T[d,d]$ is the zero vector this cannot happen; otherwise, we choose $y$ so that the degree of the restriction is as large as possible, which makes the linear inequality $3 x^\transp T[d,d] + M[d,d] \ge 0$ valid for the integer hull of $\P$, and we work inside the face of the integer hull on which this inequality is tight, where the perturbation vanishes.}



\begin{prfc}[of Proposition~\ref{prop thin ray}]
From Observation~\ref{obs f notation}, we can write $f$, \blue{up to an additive constant, which does not affect the statement,} in the form 
\[
f(x) = T[x,x,x] + M[x,x] + V[x],
\]
where $T \in \R^{n \times n \times n}$ is a symmetric tensor, $M \in \R^{n \times n}$ is a symmetric matrix, and $V \in \R^n$.
The restriction of $f$ to $\ray(y,d)$ is given by
\begin{align}
\label{eq expanded}
f_{\ray(y,d)}(\lambda) & = f(y + \lambda d) = a_3 \lambda^3 + a_2 \lambda^2 + a_1  \lambda + a_0,
\end{align}
where
\[
\begin{aligned}
a_3 & := T[d,d,d], \\
a_2 & := 3 T[y,d,d] + M[d,d], \\
a_1 & := 3 T[y,y,d] + 2 M[y,d] + V[d], \\
a_0 & := T[y,y,y] + M[y,y] + V[y].
\end{aligned}
\]

Note that $f_{\ray(y,d)}(\lambda)$ is a univariate polynomial function of degree at most three.
\blue{Consider the set of all vectors $\bar y \in \P \cap \Z^n$ such that $f$ is unbounded below on $\ray(\bar y,d)$.
This set contains $y$, and, since the degree of $f_{\ray(\bar y,d)}$ is at most three for every such $\bar y$, it contains a vector for which this degree is maximal.
It suffices to prove the statement for such a vector; thus, by replacing $y$ with it, we can assume that $y$ maximizes the degree of $f_{\ray(y,d)}$ among all vectors of $\P \cap \Z^n$ from which $f$ is unbounded below along the direction $d$, and we then show that $f$ is unbounded below on $\P \cap \Z^n \cap \pare{\ray(y,d) + \B_\epsilon}$ for every $\epsilon > 0$.}
We denote by $\deg\pare{f_{\ray(y,d)}}$ the degree of $f_{\ray(y,d)}$.
Since $f$ is not constant on $\ray(y,d)$, this degree is either 1, 2, or 3.

For every $z \in \R^n$, the restriction of $f$ to $\ray(y+z,d)$ is given by
\begin{align}
\label{eq f translated}
f(y + z + \lambda d) & = 
a_3 \lambda^3
+ \pare{a_2 + \delta_2(z)} \lambda^2 
+ \pare{a_1 + \delta_1(z)} \lambda 
+ \pare{a_0 + \delta_0(z)},
\end{align}
where
\[
\begin{aligned}
\delta_2(z) & := 3 T[z,d,d], \\
\delta_1(z) & := 6 T[y,z,d] + 3 T[z,z,d] + 2 M[z,d], \\
\delta_0(z) & := 3 T[y,z,z] + 3 T[y,y,z] + T[z,z,z] + 2 M[y,z] + M[z,z] + V[z].
\end{aligned}
\]
For every constant $\epsilon > 0$, we define the function $f^\epsilon_{\ray(y,d)} : \R_{\ge 0} \to \R$ as 
\[
f^\epsilon_{\ray(y,d)}(\lambda) := \max\bra{f(y+z+\lambda d) : \norm{z} \le \epsilon}.
\]
From \eqref{eq f translated}, for every $\epsilon > 0$ and $\lambda \ge 0$, we have
\begin{align}
\label{eq f z}
f^\epsilon_{\ray(y,d)} (\lambda) 
\le 
a_3 \lambda^3 
+ \pare{a_2 + \Delta_2} \lambda^2 + \pare{a_1 + \Delta_1} \lambda + \pare{a_0 + \Delta_0},
\end{align}
where
\[
\begin{aligned}
\Delta_2 & := 3 \, \epsilon \, \norm{T} \, \norm{d}^2, \\
\Delta_1 & := \epsilon \pare{ 6 \, \norm{T} \, \norm{y} \, \norm{d} + 2 \, \norm{M} \, \norm{d} } + 3 \, \epsilon^2 \, \norm{T} \, \norm{d}, \\
\Delta_0 & := \epsilon \pare{ 3 \, \norm{T} \, \norm{y}^2 + 2 \, \norm{M} \, \norm{y} + \norm{V} }
+ \epsilon^2 \pare{ 3 \, \norm{T} \, \norm{y} + \norm{M} }
+ \epsilon^3 \, \norm{T}.
\end{aligned}
\]
In the remainder of the proof, we consider separately the cases $\deg\pare{f_{\ray(y,d)}} \in \bra{3,2,1}$.
Recall that, by assumption, 
\begin{align}
\label{eq f limit}
\lim_{\lambda \to + \infty} f_{\ray(y,d)}(\lambda) = \lim_{\lambda \to + \infty} f(y+\lambda d) = -\infty.
\end{align}

\phead{Case $\deg\pare{f_{\ray(y,d)}} = 3$.}
In this case, $f_{\ray(y,d)}(\lambda)$ is cubic.
From \eqref{eq expanded} and \eqref{eq f limit}, we have 
\[
a_3 < 0.
\]
Fix a constant $\epsilon > 0$.
From \eqref{eq f z}, we obtain $f^\epsilon_{\ray(y,d)} (\lambda) \to -\infty$, as $\lambda \to + \infty$.
From Lemma~\ref{lem Dirichlet poly}, for every $\bar \lambda > 0$, there exists a vector in 
\[
\P \cap \Z^n \cap \bra{y + z + \lambda d : \norm{z} \le \epsilon, \lambda \ge \bar \lambda}.
\]
\blue{Since $f(y+z+\lambda d) \le f^\epsilon_{\ray(y,d)}(\lambda)$ for every $z$ with $\norm{z} \le \epsilon$, w}e obtain that $f$ is unbounded below on $\P \cap \Z^n \cap \pare{\ray(y,d) + \B_\epsilon}$.
This concludes the proof in the case $\deg\pare{f_{\ray(y,d)}} = 3$.

\phead{Case $\deg\pare{f_{\ray(y,d)}} = 2$.}
In this case, $f_{\ray(y,d)}(\lambda)$ is quadratic.
From \eqref{eq expanded} and \eqref{eq f limit}, we have 
\[
a_3 = 0, \quad a_2 < 0.
\]
Fix a constant $\epsilon > 0$.
Our goal is to show that $f$ is unbounded below on $\P \cap \Z^n \cap \pare{\ray(y,d) + \B_\epsilon}$.
\blue{Since this set shrinks as $\epsilon$ decreases, it suffices to prove the statement after reducing $\epsilon$.}
Without loss of generality, we can assume $a_2 + \Delta_2 < 0$.
To ensure this, we may reduce $\epsilon$, if $\norm{T} \neq 0$, so that
\[
\epsilon < \frac{-a_2}{3 \norm{T} \, \norm{d}^2}.
\]
\blue{(Note that when $\norm{T} = 0$ we have $\Delta_2 = 0$, and no reduction is needed.)}
From \eqref{eq f z}, we obtain $f^{\epsilon}_{\ray(y,d)} (\lambda) \to -\infty$, as $\lambda \to + \infty$.
From Lemma~\ref{lem Dirichlet poly}, for every $\bar \lambda > 0$, there exists a vector in 
\[
\P \cap \Z^n \cap \bra{y + z + \lambda d : \norm{z} \le \epsilon, \lambda \ge \bar \lambda}.
\]
We obtain that $f$ is unbounded below on $\P \cap \Z^n \cap  \pare{\ray(y,d) + \B_\epsilon}$.
This concludes the proof in the case $\deg\pare{f_{\ray(y,d)}} = 2$.

\phead{Case $\deg\pare{f_{\ray(y,d)}} = 1$.}
In this case, $f_{\ray(y,d)}(\lambda)$ is linear.
From \eqref{eq expanded} and \eqref{eq f limit}, we have 
\[
a_3 = 0, \quad a_2 = 0, \quad a_1 < 0.
\]
Fix a constant $\epsilon > 0$.
Our goal is to show that $f$ is unbounded below on $\P \cap \Z^n \cap \pare{\ray(y,d) + \B_\epsilon}$.
\blue{As in the previous case, it suffices to prove the statement after reducing $\epsilon$.}
Without loss of generality, we can assume $a_1 + \Delta_1 < 0$.
To ensure this, we may reduce $\epsilon$, so that
\[
\epsilon < 
\begin{cases}
\frac{-a_1}{\beta} & \text{if $\norm{T}=0$, $\norm{M} > 0$},\\
\frac{-\beta + \sqrt{\beta^2 - 4 a_1 \gamma}}{2 \gamma} & \text{if $\norm{T}>0$},
\end{cases}
\]
where $\beta := 6 \, \norm{T} \, \norm{y} \, \norm{d} + 2 \, \norm{M} \, \norm{d}$ and $\gamma := 3 \norm{T} \, \norm{d}$.
(Note that when $\beta = \gamma = 0$ no reduction is needed.)
We can then rewrite $a_2$ and $\delta_2(z)$ in the form
\[
\begin{aligned}
a_2 & = 3 y^\transp T[d,d] + M[d,d], \\
\delta_2(z) & = 3 z^\transp T[d,d],
\end{aligned}
\]
where we recall that $T[d,d]$ is the tensor contraction defined in Section~\ref{sec notation}.

In the remainder of the proof, we consider separately the sub-cases $T[d,d] = \z$ and $T[d,d] \neq \z$.

\phead{Sub-case $T[d,d] = \z$.}
In this sub-case we have $\delta_2(z)=0$ for every $z \in \R^n$, thus from \eqref{eq f translated}, for every $\lambda \ge 0$, we have
\[
f^\epsilon_{\ray(y,d)} (\lambda) 
\le 
a_3 \lambda^3 
+ a_2 \lambda^2 + \pare{a_1 + \Delta_1} \lambda + \pare{a_0 + \Delta_0}.
\]
Since in this case we also have $a_3 = 0$, $a_2 = 0$, and $a_1 + \Delta_1 < 0$, we obtain $f^{\epsilon}_{\ray(y,d)} (\lambda) \to -\infty$, as $\lambda \to + \infty$.
From Lemma~\ref{lem Dirichlet poly}, for every $\bar \lambda > 0$, there exists a vector in 
\[
\P \cap \Z^n \cap \bra{y + z + \lambda d : \norm{z} \le \epsilon, \, \lambda \ge \bar \lambda}.
\]
We obtain that $f$ is unbounded below on $\P \cap \Z^n \cap \pare{\ray(y,d) + \B_\epsilon}$.
This concludes the proof in this sub-case.

\phead{Sub-case $T[d,d] \neq \z$.}
In this sub-case, we show that, for every $x \in \P \cap \Z^n$, we have 
\[
3 x^\transp T[d,d] + M[d,d] \ge 0.
\]
Assume, for a contradiction, that there exists $x \in \P \cap \Z^n$ such that $3 x^\transp T[d,d] + M[d,d] < 0$.
Then, $\ray(x,d)$ is a ray of $\P$.
The restriction of $f$ to $\ray(x,d)$ is
\[
f_{\ray(x,d)}(\lambda) = f(x + \lambda d) = a_3 \lambda^3 + a'_2 \lambda^2 + a'_1 \lambda + a'_0,
\]
where
\[
\begin{aligned}
a'_2 & := 3 T[x,d,d] + M[d,d], \\
a'_1 & := 3 T[x,x,d] + 2 M[x,d] + V[d], \\
a'_0 & := T[x,x,x] + M[x,x] + V[x].
\end{aligned}
\]
Since $a_3=0$ and $a'_2 < 0$, $f$ is unbounded below on $\ray(x,d)$.
However, \blue{$x \in \P \cap \Z^n$ and} $\deg\pare{f_{\ray(x,d)}} = 2 > 1 = \deg\pare{f_{\ray(y,d)}}$, which contradicts our choice of $y$.

Let $\P_I$ be the integer hull of $\P$, i.e., $\P_I = \conv(\P \cap \Z^n)$.
It is well known that $\P_I$ is a rational polyhedron\blue{; see theorem~16.1 in~\cite{SchBookIP}}.
The linear inequality $3 x^\transp T[d,d] + M[d,d] \ge 0$ is then also valid for $\P_I$.
Consider now the hyperplane
\[
\begin{aligned}
\H :& = \bra{x \in \R^n : 3 x^\transp T[d,d] + M[d,d] = 0} \\
& = \bra{y+z \in \R^n : z^\transp T[d,d] = 0},
\end{aligned}
\]
and the face $\F$ of $\P_I$ defined by 
\[
\F := \P_I \cap \H.
\]
\blue{Note that $\F$ is a rational polyhedron, since it is a face of the rational polyhedron $\P_I$, even though $\H$ need not be a rational hyperplane; moreover, $\F \cap \Z^n \subseteq \P_I \cap \Z^n = \P \cap \Z^n$.}

We show that the ray $\ray(y,d)$ is contained in $\F$. 
To do so, observe that $\ray(y,d)$ is contained in $\P_I$, because $y \in \P_I$ and $d \in \rec(\P) = \rec(\P_I)$.
Furthermore, $\ray(y,d)$ is contained in $\H$ because
\[
\begin{aligned}
3 T[y+\lambda d,d,d] + M[d,d]
& = 3 T[y,d,d] + 3 \lambda T[d,d,d] + M[d,d] \\
& = 3 \lambda a_3 + a_2 \\
& = 0,
\end{aligned}
\]
where the last equality follows from the assumptions $a_3=a_2=0$ of this case.

Define the function $f^{\epsilon,\F}_{\ray(y,d)} : \R_{\ge 0} \to \R$ as 
\[
f^{\epsilon,\F}_{\ray(y,d)}(\lambda) := \max\bra{f(y+z+\lambda d) : z^\transp T[d,d] = 0 , \, \norm{z} \le \epsilon}.
\]
From \eqref{eq f translated}, for every $\lambda \ge 0$, we have
\[
f^{\epsilon,\F}_{\ray(y,d)} (\lambda)
\le 
a_3 \lambda^3 
+ a_2 \lambda^2 + \pare{a_1 + \Delta_1} \lambda + \pare{a_0 + \Delta_0}.
\]
Since in this case we also have $a_3 = 0$, $a_2 = 0$, and $a_1 + \Delta_1 < 0$, we obtain $f^{\epsilon,\F}_{\ray(y,d)} (\lambda) \to -\infty$, as $\lambda \to + \infty$.
From Lemma~\ref{lem Dirichlet poly}, applied to $\F$, for every $\bar \lambda > 0$, there exists a vector in 
\[
\begin{aligned}
& \F \cap \Z^n \cap \bra{y + z + \lambda d : \norm{z} \le \epsilon, \, \lambda \ge \bar \lambda} \\
& = \F \cap \Z^n \cap \bra{y + z + \lambda d : z^\transp T[d,d] = 0 , \, \norm{z} \le \epsilon, \, \lambda \ge \bar \lambda}.
\end{aligned}
\]
\blue{Here the equality of the two sets holds because every vector of $\F \subseteq \H$ at distance less than $\epsilon$ from a vector $y + \lambda d \in \H$ can be written as $y + \lambda d + z$ with $z^\transp T[d,d] = 0$ and $\norm{z} < \epsilon$.}
We obtain that $f$ is unbounded below on $\F \cap \Z^n \cap \pare{\ray(y,d) + \B_\epsilon} \subseteq \P \cap \Z^n \cap \pare{\ray(y,d) + \B_\epsilon}$.
This concludes the proof in this sub-case.
\end{prfc}

\section{Proofs of Propositions~\blue{\ref{prop quartic example}}, \ref{prop cubic example}, and~\ref{prop quadratic example}}
\label{sec other}

\subsection{Proof of Proposition~\ref{prop quartic example}}
\label{sec quartic example}

\begin{bluepart}
Let $k \ge 4$ be an integer.
We show that Proposition~\ref{prop quartic example} holds with
\begin{align}
\label{eq quartic example}
\P := \bra{x \in \R^2 : x_1 \le 0}
\qquad \text{and} \qquad
f(x) := x_1 + \pare{1 + (-x_1)^{k-4}} \pare{x_1^2 - x_2}^2,
\end{align}
with the convention $(-x_1)^0 = 1$.
For $k = 4$, up to a factor two in the second term, $f$ is the function $x_1 + (x_1^2 - x_2)^2$, which tends to $-\infty$ along the integer points $(-t,t^2)$, $t \in \N$, of the parabola $x_2 = x_1^2$, while the region where it is nonpositive is a neighborhood of this parabola whose width grows only like $\sqrt{\abs{x_1}}$, so that it eventually leaves every strip of bounded width.
The factor $1 + (-x_1)^{k-4}$, which is at least one on $\P$, only raises the degree to $k$ without affecting this picture.
The function $x_1 + (x_1^2 - x_2)^2$ is a variant of the quartic $y_2 - (y_2 - y_1^2)^2$ in proposition~5.3 in~\cite{BiedPHil23MPB}, which is unbounded on $\R^2$ but bounded on every ray, and which is in turn inspired by a classical example of Frank and Wolfe~\cite{FraWol56}; Proposition~\ref{prop quartic example} adds the integer setting, the extension to every degree $k \ge 4$, and the boundedness on neighborhoods of rays of any width.
We first prove a slightly more general statement, over the reals rather than only on the integer points, which is also what the comparison with the continuous setting in Section~\ref{sec intro} uses.

\begin{lemma}
\label{lem parabola}
Let $\P \subseteq \R^2$ be a polyhedron, let $h \colon \R^2 \to \R$ be a polynomial function such that $h(x) \ge 1$ for every $x \in \P$, and let
\[
f(x) := x_1 + h(x) \pare{x_1^2 - x_2}^2.
\]
Then, for every ray $\ray(y,d)$ of $\P$ and every constant $\epsilon > 0$, $f$ is bounded below on $\P \cap \pare{\ray(y,d) + \B_\epsilon}$.
\end{lemma}

\begin{prf}
Fix a ray $\ray(y,d)$ of $\P$ and $\epsilon > 0$, and let $\mathcal S := \bra{x \in \P : f(x) \le 0}$.
We first show that $\mathcal S$ is contained in a thin region around the parabola $x_2 = x_1^2$, $x_1 \le 0$.
Let $x \in \mathcal S$.
Since $h(x) \ge 1$, we have $x_1 \le - h(x) \pare{x_1^2 - x_2}^2 \le 0$, and $\pare{x_1^2 - x_2}^2 \le - x_1 / h(x) \le - x_1$.
Writing $x_1 = -t$ with $t \ge 0$, we obtain
\begin{align}
\label{eq parabola region}
x_2 = t^2 + s \qquad \text{for some } s \in \R \text{ with } \abs{s} \le \sqrt t.
\end{align}

Next, we show that $\mathcal S \cap \pare{\ray(y,d) + \B_\epsilon}$ is bounded.
Let $\mathcal L := \bra{y + \lambda d : \lambda \in \R}$ be the line containing $\ray(y,d)$, so that $\ray(y,d) + \B_\epsilon \subseteq \mathcal L + \B_\epsilon$, and it suffices to show that $\mathcal S \cap \pare{\mathcal L + \B_\epsilon}$ is bounded.
We distinguish two cases.

Assume first $d_1 = 0$, so that $\mathcal L = \bra{x \in \R^2 : x_1 = y_1}$ is a vertical line.
For every $x \in \mathcal S \cap \pare{\mathcal L + \B_\epsilon}$, we have $\abs{x_1 - y_1} \le \epsilon$, thus $t = -x_1 \le \abs{y_1} + \epsilon =: t_0$, and, by \eqref{eq parabola region}, $\abs{x_2} \le t_0^2 + \sqrt{t_0}$.
Hence $\mathcal S \cap \pare{\mathcal L + \B_\epsilon}$ is bounded.

Assume now $d_1 \neq 0$.
Then $\mathcal L = \bra{x \in \R^2 : x_2 = a x_1 + b}$, where $a := d_2/d_1$ and $b := y_2 - a y_1$.
The Euclidean distance of a point $x \in \R^2$ from $\mathcal L$ is $\abs{x_2 - a x_1 - b} / \sqrt{1 + a^2}$, thus every $x \in \mathcal L + \B_\epsilon$ satisfies
\[
\abs{x_2 - a x_1 - b} \le \epsilon' := \epsilon \sqrt{1+a^2}.
\]
On the other hand, for every $x \in \mathcal S$, writing $x_1 = -t$ and using \eqref{eq parabola region},
\[
x_2 - a x_1 - b = t^2 + s + a t - b \ \ge\ t^2 - \sqrt t - \abs{a} t - \abs{b} =: \phi(t).
\]
Since $\phi(t) \to +\infty$ as $t \to +\infty$, there exists $t_0 \ge 0$ such that $\phi(t) > \epsilon'$ for every $t > t_0$.
Therefore, every $x \in \mathcal S \cap \pare{\mathcal L + \B_\epsilon}$ satisfies $t = -x_1 \le t_0$, and, by \eqref{eq parabola region}, $\abs{x_2} \le t_0^2 + \sqrt{t_0}$.
Hence $\mathcal S \cap \pare{\mathcal L + \B_\epsilon}$ is bounded.

We can now conclude.
Since $\mathcal S \cap \pare{\ray(y,d) + \B_\epsilon}$ is bounded, \blue{its closure is compact, and since $f$ is continuous,} $f$ is bounded below on this set by some constant $c \in \R$.
Every $x \in \P \cap \pare{\ray(y,d) + \B_\epsilon}$ that is not in $\mathcal S$ satisfies $f(x) > 0$.
Therefore, $f \ge \min\bra{c,0}$ on $\P \cap \pare{\ray(y,d) + \B_\epsilon}$.
\end{prf}

We are now ready to prove Proposition~\ref{prop quartic example}.

\begin{prfc}[of Proposition~\ref{prop quartic example}]
Let $\P$ and $f$ be as in \eqref{eq quartic example}.
Note that $\P$ is a rational polyhedron.
The function $f$ is a polynomial with integer, hence rational, coefficients.
Its degree is exactly $k$: for $k=4$, we have $f(x) = x_1 + 2 \pare{x_1^2 - x_2}^2$, which has degree four, while for $k \ge 5$ the monomials of $(-x_1)^{k-4} \pare{x_1^2 - x_2}^2$ have degrees $k$, $k-1$, and $k-2$, the one of degree $k$ being $(-1)^{k-4} x_1^{k}$, and all other monomials of $f$ have degree at most four.

We first show that $f$ is unbounded below on $\P \cap \Z^2$.
For every $t \in \N$, the vector $x^t := (-t, t^2)$ is in $\P \cap \Z^2$, and $(x^t_1)^2 - x^t_2 = 0$, thus $f(x^t) = -t$.
Hence $f(x^t) \to -\infty$ as $t \to +\infty$.

Next, let $h(x) := 1 + (-x_1)^{k-4}$.
For every $x \in \P$, we have $-x_1 \ge 0$, thus $h(x) \ge 1$.
By Lemma~\ref{lem parabola}, for every ray $\ray(y,d)$ of $\P$ and every $\epsilon > 0$, $f$ is bounded below on $\P \cap \pare{\ray(y,d) + \B_\epsilon}$.
\end{prfc}

We remark that the inequality $x_1 \le 0$ in the definition of $\P$ is needed only when $k$ is odd, since $(-x_1)^{k-4}$ is then negative for $x_1 > 0$.
If $k$ is even, we have $(-x_1)^{k-4} = x_1^{k-4} \ge 0$ and $h \ge 1$ on all of $\R^2$, thus Lemma~\ref{lem parabola} applies with $\P = \R^2$, and Proposition~\ref{prop quartic example} holds with $\P = \R^2$; in particular, for $k=4$ it holds with $\P = \R^2$ and $f(x) = x_1 + \pare{x_1^2 - x_2}^2$.
\end{bluepart}

\subsection{Proof of Proposition~\ref{prop cubic example}}
\label{sec cubic example}


\blue{The idea is to take a cubic form that is nonnegative on $\rec(\P)$ and vanishes only on one ray of $\rec(\P)$, whose direction $d$ is irrational, together with a linear term that is negative in the direction $d$.
Along every ray of $\P$ whose direction is not a positive multiple of $d$, the cubic term dominates and $f$ tends to $+\infty$, while along the ray of $\P$ that starts at the origin with direction $d$ the cubic term vanishes and the linear term drives $f$ to $-\infty$; Proposition~\ref{prop thin ray} then yields the unboundedness of $f$ on $\P \cap \Z^3$.
A ray with irrational direction contains at most one integer point, so $f$ is bounded below also on the integer points of the rays with direction $d$, and therefore on the integer points of every ray of $\P$.
The construction is based on the one in proposition~2.11 in~\cite{BiedPHil23MPB}, which is a statement over the reals: we use the same homogeneous cubic, up to sign, and the same irrational direction $\tilde d := (2^{1/3},2^{2/3},1)$, but a different box, a linear rather than quadratic perturbation, and a proof, by the AM--GM inequality, that $T[d,d,d] \ge 0$ on $\Qpoly$ with equality only at $\tilde d$, in place of example~2.1 and observation~2.2 there.
Proposition~\ref{prop cubic example} is tight with respect to the dimension: for $n \le 2$, theorem~1.3 in~\cite{dPHilWeiZem16MOR} certifies unboundedness by exhibiting an integer point of $\P$ and an integer vector of $\rec(\P)$ along which $f$ is unbounded below, so that rays suffice.}

\begin{prfc}[\blue{of Proposition~\ref{prop cubic example}}]
Let $T \in \Z^{3 \times 3 \times 3}$ be the symmetric tensor with nonzero entries
\[
\begin{aligned}
& T_{111}=2, \ T_{222}=1, \ T_{333}=4, \\
& T_{123} = T_{132} = T_{213} = T_{231} = T_{312} = T_{321} = -1, \\
\end{aligned}
\]
and let $V \in \Z^{3}$ with nonzero entry $V_1 = -1$.
Let $f : \R^3 \to \R$ be the cubic function defined by 
\[
f(x) := T[x,x,x] + V[x].
\]

Let $\P \subseteq \R^3$ be the rational polyhedron
\[
\P := \cone \Qpoly, \quad \text{where} \quad
\Qpoly := \{x \in \R^3 : 1 \le x_1 \le 2, \ 1 \le x_2 \le 2, \ x_3 = 1 \}.
\]

Let $\ray(y,d)$ be a ray of $\P$\blue{, so that $y \in \P$ and $d \in \rec(\P) \setminus \bra{\z} = \P \setminus \bra{\z}$.
Since every nonzero vector of $\cone \Qpoly$ is a positive multiple of a vector of $\Qpoly$, and scaling $d$ by a positive constant does not change $\ray(y,d)$, we can assume without loss of generality that $d \in \Qpoly$}.
The restriction of $f$ to $\ray(y,d)$ is 
\[
\begin{aligned}
f_{\ray(y,d)}(\lambda) & = f(y + \lambda d) \\
& = T[d,d,d] \lambda^3 + 3 T[y,d,d] \lambda^2 + \pare{3 T[y,y,d] + V[d]}  \lambda + f(y).
\end{aligned}
\]


We show that $T[d,d,d] \ge 0$ for every $d \in \Qpoly$.
Since $d_3 = 1$ on $\Qpoly$, we have
\[
T[d,d,d] = 2d_1^3 + d_2^3 + 4 - 6d_1 d_2.
\]
By the AM--GM inequality applied to the three positive terms $2d_1^3$, $d_2^3$, and $4$, we obtain
\[
2d_1^3 + d_2^3 + 4 \ge 3 \bigl(2d_1^3 \cdot d_2^3 \cdot 4\bigr)^{1/3} = 3 \bigl(8 d_1^3 d_2^3\bigr)^{1/3} = 6 d_1 d_2,
\]
with equality if and only if $2d_1^3 = d_2^3 = 4$, i.e., $d_1 = 2^{1/3}$ and $d_2 = 2^{2/3}$.
Therefore, $T[d,d,d] \ge 0$ for every $d \in \Qpoly$, and $T[d,d,d] = 0$ holds only when $d$ is the irrational vector $\tilde d := (2^{1/3}, 2^{2/3}, 1)$.

Hence, if $d \neq \tilde d$, we have $f_{\ray(y,d)}(\lambda) \to + \infty$, as $\lambda \to + \infty$.
In particular, if $d \neq \tilde d$, $f$ is bounded below on $\ray(y, d)$, and so also on $\ray(y, d) \cap \Z^3$.

Furthermore, $V[\tilde d] = - \tilde d_1 = -2^{1/3} < 0$, thus $f_{\ray(\z,\tilde d)}(\lambda) \to - \infty$, as $\lambda \to + \infty$.
It then follows from Proposition~\ref{prop thin ray} that $f$ is unbounded below on $\P \cap \Z^3$.

It remains to show that $f$ is bounded below on $\ray(y, \tilde d) \blue{\cap \Z^3}$ for every $y \in \P$. 
This follows from the irrationality of $\tilde d$, which implies that each such ray contains at most one integer point\blue{: two of them would give $(\lambda_1 - \lambda_2) \tilde d \in \Z^3$ for some $\lambda_1 \neq \lambda_2$, and since $\tilde d_3 = 1$ this would force $\tilde d \in \Q^3$}.
\end{prfc}

\subsection{Proof of Proposition~\ref{prop quadratic example}}

\blue{The construction follows the same template as the one in the proof of Proposition~\ref{prop cubic example}, with a quadratic form in place of the cubic one; the argument is simpler, because the nonnegativity of the leading coefficient on $\Qpoly$ now follows by completing the square.}

\begin{prfc}[\blue{of Proposition~\ref{prop quadratic example}}]
Let $M \in \R^{2 \times 2}$ be the symmetric matrix
\[
M = 
\begin{pmatrix}
1 & -\sqrt{2} \\
-\sqrt{2} & 2
\end{pmatrix}
\]
and let $V \in \R^2$ with nonzero entry $V_1 = -1$.
Let $f : \R^2 \to \R$ be the quadratic function defined by 
\[
f(x) := M[x,x] + V[x].
\]

Let $\P \subseteq \R^2$ be the rational polyhedron
\[
\P := \cone \Qpoly, \quad \text{where} \quad
\Qpoly := \{x \in \R^2 : 1 \le x_1 \le 2, \ x_2 = 1 \}.
\]

Let $\ray(y,d)$ be a ray of $\P$\blue{, so that $y \in \P$ and $d \in \rec(\P) \setminus \bra{\z} = \P \setminus \bra{\z}$.
Since every nonzero vector of $\cone \Qpoly$ is a positive multiple of a vector of $\Qpoly$, and scaling $d$ by a positive constant does not change $\ray(y,d)$, we can assume without loss of generality that $d \in \Qpoly$}.
The restriction of $f$ to $\ray(y,d)$ is 
\[
\begin{aligned}
f_{\ray(y,d)}(\lambda) & = f(y + \lambda d) \\
& = M[d,d] \lambda^2 + \pare{2 M[y,d] + V[d]}  \lambda + f(y).
\end{aligned}
\]


We show that $M[d,d] \ge 0$ for every $d \in \Qpoly$.
Since $d_2 = 1$ on $\Qpoly$, we have
\[
M[d,d] = d_1^2 - 2\sqrt{2}\, d_1 + 2 = \bigl(d_1 - \sqrt{2}\bigr)^2 \ge 0,
\]
with equality if and only if $d_1 = \sqrt{2}$.
Therefore, $M[d,d] \ge 0$ for every $d \in \Qpoly$, and $M[d,d] = 0$ holds only when $d$ is the irrational vector $\tilde d := (\sqrt{2}, 1)$.

Hence, if $d \neq \tilde d$, we have $f_{\ray(y,d)}(\lambda) \to + \infty$, as $\lambda \to + \infty$.
In particular, if $d \neq \tilde d$, $f$ is bounded below on $\ray(y, d)$, and so also on $\ray(y, d) \cap \Z^2$.

Furthermore, $V[\tilde d] = - \tilde d_1 = -\sqrt{2} < 0$, thus $f_{\ray(\z,\tilde d)}(\lambda) \to - \infty$, as $\lambda \to + \infty$.
It then follows from Proposition~\ref{prop thin ray} that $f$ is unbounded below on $\P \cap \Z^2$.

It remains to show that $f$ is bounded below on $\ray(y, \tilde d) \blue{\cap \Z^2}$ for every $y \in \P$. 
This follows from the irrationality of $\tilde d$, which implies that each such ray contains at most one integer point\blue{: two of them would give $(\lambda_1 - \lambda_2) \tilde d \in \Z^2$ for some $\lambda_1 \neq \lambda_2$, and since $\tilde d_2 = 1$ this would force $\tilde d \in \Q^2$}.
\end{prfc}

\section{\blue{Proof of Corollary~\ref{th quad}}}
\label{sec proof cor quad}

\blue{Corollary~\ref{th quad} is the special case of Theorem~\ref{th main} in which $f$ is a quadratic function.
We give a direct proof of it that does not rely on Proposition~\ref{prop int ray}, and that only uses Proposition~\ref{prop thin ray} together with the continuous theory of unbounded quadratic optimization; besides being of independent interest, it serves as a warm-up for the proof of Proposition~\ref{prop int ray} in Section~\ref{proof prop int ray}.
The proof consists of two steps.
In the first step, we obtain a ray $\ray(y,d)$ of $\P$, with $y$ not necessarily integer, along which $f$ is unbounded below; in the quadratic case this ray is provided by theorem~5.2 in~\cite{BiedPHil23MPB}, which holds for real coefficients.
In the second step, we transfer the unboundedness to a ray $\ray(\bar z,d)$ starting at an integer point $\bar z \in \P \cap \Z^n$; the tool for this step is the integrality of $\P$, which allows us to find an integer point of $\P$ on which an affine function that is negative somewhere on $\P$ is negative as well.
The overall structure follows the one of proposition~4 in~\cite{dPDeyMol17MPA}, with two changes: the real ray comes from theorem~5.2 in~\cite{BiedPHil23MPB} rather than from~\cite{Vav90}, so that rationality is not needed, and Proposition~\ref{prop thin ray} replaces the final step there, which relies on the rationality of the direction of the ray.}

\begin{prfc}[\blue{of Corollary~\ref{th quad}}]
\begin{bluepart}
If there exists a ray $\ray(y,d)$ of $\P$ such that, for every $\epsilon > 0$, $f$ is unbounded below on $\P \cap \Z^n \cap \pare{\ray(y,d) + \B_\epsilon}$, then $f$ is unbounded below on $\P \cap \Z^n$, since $\P \cap \Z^n \cap \pare{\ray(y,d) + \B_\epsilon} \subseteq \P \cap \Z^n$.
This proves the ``if'' direction.

For the ``only if'' direction, assume that $f$ is unbounded below on $\P \cap \Z^n$.
Due to Proposition~\ref{prop thin ray}, it suffices to show that there exists a ray $\ray(\bar z,d)$ of $\P$ with $\bar z \in \Z^n$ such that $f$ is unbounded below on $\ray(\bar z,d)$.

We first show that we can assume that $\P$ is integral, i.e., that $\P = \conv(\P \cap \Z^n)$.
Let $\P_I := \conv(\P \cap \Z^n)$ be the integer hull of $\P$.
Since $\P$ is a rational polyhedron and $\P \cap \Z^n$ is nonempty, because $f$ is unbounded below on it, it is well known that $\P_I$ is a rational polyhedron with $\P_I \cap \Z^n = \P \cap \Z^n$ and $\rec(\P_I) = \rec(\P)$; see theorem~16.1 in~\cite{SchBookIP}.
In particular, $f$ is unbounded below on $\P_I \cap \Z^n$, and every ray of $\P_I$ is a ray of $\P$, since $\P_I \subseteq \P$ and $\rec(\P_I) = \rec(\P)$.
Hence it suffices to prove the statement for $\P_I$, and we assume from now on that $\P$ is integral.

From Observation~\ref{obs f notation}, we can write $f$, up to an additive constant, which does not affect the statement, in the form
\[
f(x) = M[x,x] + V[x],
\]
where $M \in \R^{n \times n}$ is a symmetric matrix and $V \in \R^n$.
Recall from Section~\ref{sec notation} that $M[d] \in \R^n$ denotes the contraction of $M$ with $d$, so that $M[z,d] = z^\transp M[d]$ for every $z \in \R^n$.

Since $\P \cap \Z^n \subseteq \P$, the function $f$ is unbounded below on $\P$.
It then follows from theorem~5.2 in~\cite{BiedPHil23MPB} that there exists a ray $\ray(y,d)$ of $\P$ such that $f$ is unbounded below on $\ray(y,d)$.
Here $y \in \P$ and $d \in \rec(\P) \setminus \bra{\z}$, but $y$ need not be an integer vector.
For every $z \in \R^n$, by expanding the bilinear form and the linear form, the restriction of $f$ to $\ray(z,d)$ is given by
\begin{align}
\label{eq quad restriction}
\begin{aligned}
f_{\ray(z,d)}(\lambda)
&= M[z + \lambda d,\, z + \lambda d] + V[z + \lambda d] \\
&= M[z,z] + 2\lambda\, M[z,d] + \lambda^2 M[d,d] + V[z] + \lambda V[d] \\
&= M[d,d]\,\lambda^2 + g(z) \, \lambda + f(z),
\qquad \text{where} \quad g(z) := 2\, z^\transp M[d] + V[d],
\end{aligned}
\end{align}
and we note that $g : \R^n \to \R$ is an affine function of $z$.
Since $f_{\ray(y,d)}$ is a univariate polynomial of degree at most two that is unbounded below on $\R_{\ge 0}$, exactly one of the following two cases occurs: either $M[d,d] < 0$, or $M[d,d] = 0$ and $g(y) < 0$.
Indeed, if $M[d,d] > 0$, or if $M[d,d] = 0$ and $g(y) \ge 0$, then $f_{\ray(y,d)}$ is bounded below on $\R_{\ge 0}$.

We now transfer this unboundedness to a ray of $\P$ that starts at an integer point.
For every $\bar z \in \P \cap \Z^n$, the set $\ray(\bar z,d)$ is a ray of $\P$, because $\bar z \in \P$ and $d \in \rec(\P)$, and its restriction is given by \eqref{eq quad restriction} with $z = \bar z$.
Note that the quadratic coefficient $M[d,d]$ does not depend on the starting point of the ray, while the linear coefficient $g(\bar z)$ does.

Consider first the case $M[d,d] < 0$.
Then $f_{\ray(\bar z,d)}$ is a quadratic polynomial with negative leading coefficient, thus $f_{\ray(\bar z,d)}(\lambda) \to -\infty$ as $\lambda \to +\infty$, for every $\bar z \in \P \cap \Z^n$, and we are done.

Consider now the case $M[d,d] = 0$ and $g(y) < 0$.
In this case, $f_{\ray(z,d)}(\lambda) = g(z) \lambda + f(z)$ for every $z \in \R^n$, thus it suffices to find $\bar z \in \P \cap \Z^n$ with $g(\bar z) < 0$.
Since $\P$ is integral, we have $y \in \P = \conv(\P \cap \Z^n)$, thus there exist $z_1,\dots,z_r \in \P \cap \Z^n$ and $\mu_1,\dots,\mu_r \ge 0$ with $\sum_{l=1}^r \mu_l = 1$ such that $y = \sum_{l=1}^r \mu_l z_l$.
Since $g$ is affine, we obtain
\[
\sum_{l=1}^r \mu_l \, g(z_l) = g \pare{\sum_{l=1}^r \mu_l z_l} = g(y) < 0,
\]
hence $g(z_l) < 0$ for at least one $l \in \bra{1,\dots,r}$.
Let $\bar z := z_l$.
Then $f_{\ray(\bar z,d)}(\lambda) = g(\bar z) \lambda + f(\bar z) \to -\infty$ as $\lambda \to +\infty$, and we are done.
(We remark that, if the contraction $M[d]$ is the zero vector, then $g$ is the constant $V[d] < 0$, and any $\bar z \in \P \cap \Z^n$ works.)
\end{bluepart}

\blue{In both cases, we have found a ray $\ray(\bar z,d)$ of $\P$ with $\bar z \in \Z^n$ such that $f$ is unbounded below on $\ray(\bar z,d)$.
Proposition~\ref{prop thin ray} then implies that there exists $\bar y \in \P \cap \Z^n$ such that, for every $\epsilon > 0$, $f$ is unbounded below on $\P \cap \Z^n \cap \pare{\ray(\bar y,d) + \B_\epsilon}$, which concludes the proof.}
\end{prfc}

\blue{This argument does not extend to cubic functions, although theorem~5.2 in~\cite{BiedPHil23MPB} holds for polynomials of degree at most three, so that a ray of $\P$ along which $f$ is unbounded below is still available.
The transfer of the unboundedness to an integer starting point still works when the restriction of $f$ to the ray is cubic or quadratic, because the corresponding leading coefficient does not depend on the starting point, or depends on it affinely.
It breaks down when the restriction is linear: there the coefficient of $\lambda$ is quadratic in the starting point, so writing $y$ as a convex combination of integer points of $\P$ no longer yields one at which that coefficient is negative.
For this reason, the proof of Proposition~\ref{prop int ray} in Section~\ref{proof prop int ray} does not start from a ray over the reals, and it obtains a direction as a limit of normalized integer points instead.}

\section{Proof of Proposition~\ref{prop int ray}}
\label{proof prop int ray}

\blue{In this section, we prove Proposition~\ref{prop int ray}, thereby completing the proof of Theorem~\ref{th main}.
The proof is quite long.
To make it more readable, it is subdivided into several claims.}

\begin{greenpart}
The proof is by induction on $n$, and the inductive step consists of three stages.
In the first stage, we reduce to the case in which $\P$ is pointed, integral, and full-dimensional, and contains the origin in its interior (Claims~\ref{claim pointed}--\ref{claim origin}), and we extract a sequence of integer points of $\rec(\P)$ along which $f$ tends to $-\infty$ (Claim~\ref{claim sequence}).
At this point, either $f$ is unbounded below on some ray of $\P$ starting at an integer point, and the proposition holds, or $f$ is bounded below on every such ray; the remaining two stages are carried out in the latter case, and show that it cannot occur.
In the second stage, we establish some elementary facts about the restrictions of $f$ to such rays (Claims~\ref{claim sign}--\ref{claim C}).
In the third stage, we iterate a projection argument: we successively project the sequence, along directions obtained as limits of normalized points, onto faces of $\rec(\P)$ of smaller and smaller dimension, and we show that the projected sequences inherit the unboundedness of $f$.
The projection can always be iterated (Claim~\ref{claim F}), while the faces that it produces always have dimension at least two (Claim~\ref{claim dim}); since each projection strictly decreases the dimension, these two facts are incompatible, and the proposition follows.
The projected points are, in general, not integer.
This is not an obstacle, because the argument uses the integrality of the original sequence and of $\P$, and never that of the projected points; this is what allows the projection to be iterated.

The sequence provided by Claim~\ref{claim sequence} is not the sequence used in the continuous setting.
There, the corresponding statement is lemma~5.1 in~\cite{BiedPHil23MPB}, which is due to~\cite{AndBelShi82,Kla19}, and which yields a sequence $\bra{x^j}_{j \in \N}$ of vectors of $\P$ along which $f$ decreases to $-\infty$, whose norms increase to $+\infty$, and which is minimal, in the sense that $f(x) > f(x^j)$ for every $x \in \P$ with $\norm{x} < \norm{x^j}$.
We do not require the minimality property, and we require instead that the vectors be integer and lie in $\rec(\P)$.
These two properties have no counterpart in the continuous setting, and they are what makes the sequence a chain of order~$0$, in the sense defined below, which is the object that the projection argument of the third stage operates on.
\end{greenpart}

\begin{prfc}[\blue{of Proposition~\ref{prop int ray}}]
Let $\P = \bra{x \in \R^n : Ax \le b}$, where $A \in \Z^{m \times n}$ and $b \in \Z^m$\green{; this is without loss of generality, since $\P$ is a rational polyhedron, and we can scale each of its defining inequalities by a positive integer}.
We then have $\rec(\P) = \bra{x \in \R^n : Ax \le \z}$.

The proof is by induction on $n$.
If $n=1$, then $\P$ is unbounded, since $f$ is unbounded below on $\P \cap \Z$, thus either $\P = \R$, or $\P = \bra{ y + \lambda d : \lambda \ge 0}$ with $y \in \Q$ and $d \in \bra{+1,-1}$.
\green{In the first case, the univariate polynomial $f$ is unbounded below on $\Z$, thus it tends to $-\infty$ as $\lambda \to +\infty$ or as $\lambda \to -\infty$, and the proposition holds by choosing the ray $\ray(0,1)$ or the ray $\ray(0,-1)$, respectively.}
In the second case, \green{the univariate polynomial $\lambda \mapsto f(y + \lambda d)$ is unbounded below on $[0,+\infty)$, thus it tends to $-\infty$ as $\lambda \to + \infty$, and} the proposition holds, by choosing the ray $\ray(\ceil{y},1)$, if $d = +1$, and the ray $\ray(\floor{y},-1)$, if $d = -1$.
This concludes our base case. In the remainder of the proof, we consider the inductive case.

\phead{\blue{Stage 1: reductions.}}
\begin{claim}
\label{claim pointed}
We can assume that $\P$ is pointed.
\end{claim}

\begin{cpf}
Let $\P_1,\P_2,\dots,\P_{2^n}$ be the intersections of $\P$ with the $2^n$ orthants of $\R^n$.
Since each vector in $\P \cap \Z^n$ is in at least one of these finitely many polyhedra, there exists $i \in \bra{1,2,\dots,2^n}$ such that $f$ is unbounded below on $\P_i \cap \Z^n$.
\green{Each $\P_i$ is pointed, since $\rec(\P_i)$ is contained in an orthant, and therefore contains no line.}
Furthermore, $\rec(\P_i) \subseteq \rec(\P)$, for every $i \in \bra{1,2,\dots,2^n}$\blue{, thus every ray of $\P_i$ is a ray of $\P$}.
\end{cpf}

\begin{claim}
\label{claim integral}
We can assume that $\P$ is integral.
\end{claim}

\begin{cpf}
It suffices to recall that the integer hull of $\P$, which is the convex hull of the integer points in $\P$, is a rational polyhedron with the same recession cone as $\P$\blue{; see theorem~16.1 in~\cite{SchBookIP}.
Moreover, it has the same integer points as $\P$, and every ray of it is a ray of $\P$}.
\end{cpf}

\begin{claim}
\label{claim H}
We can assume that, for every rational hyperplane $\H$ of $\R^n$, $f$ is bounded below on $\P \cap \H \cap \Z^n$.
In particular, we can assume that $\P$ is full-dimensional.
\end{claim}

\begin{cpf}
Assume that there exists a rational hyperplane $\H$ of $\R^n$ such that $f$ is unbounded below on $\P \cap \H \cap \Z^n$.
In particular, $\H$ contains integer points.
It then follows from Lemma~\ref{lem projection} that there exists a map $\pi : \R^{n-1} \to \R^n$ of the form $\pi(y) = \tilde x + \blue{U} y$ with $\tilde x \in \Z^n$ and $\blue{U} \in \Q^{n \times (n-1)}$ of full rank, such that
\[
\begin{aligned}
\H & = \pi(\R^{n-1}), \\
\H \cap \Z^n & = \pi(\Z^{n-1}).
\end{aligned}
\]
Let
\[
\P' :=  \bra{x' \in \R^{n-1} : A {\bluecolor U} x' \le b - A \tilde x},
\]
so that $\pi(\P') = \P \cap \H$.
Furthermore, let $f' : \R^{n-1} \to \R$ be defined by $f'(x') := f(\pi(x')) = f(\tilde x + \blue{U} x')$.
Then, $f'$ is unbounded below on $\P' \cap \Z^{n-1}$\green{, since $U$ has full rank, thus $\pi$ is injective, and $\pi(\P' \cap \Z^{n-1}) = \pi(\P') \cap \pi(\Z^{n-1}) = \P \cap \H \cap \Z^n$}.
By induction, there exists a ray $\ray(y',d')$ of $\P'$ with $y' \in \Z^{n-1}$ such that $f'$ is unbounded below on $\ray(y',d')$.
Let $y:=\tilde x + \blue{U} y'$ and $d := \blue{U}d'$.
Then, $y \in \Z^n$ and $\ray(y,d)$ is a ray of $\P$.
Furthermore, $f$ is unbounded below on this ray, so we are done.
Therefore, we can assume that, for every rational hyperplane $\H$ of $\R^n$, $f$ is bounded below on $\P \cap \H \cap \Z^n$.

If $\P$ is not full-dimensional, then it is contained in a rational hyperplane $\H$.
But then $f$ is bounded below on $\P \cap \H \cap \Z^n = \P \cap \Z^n$, a contradiction.
\end{cpf}

\begin{claim}
\label{claim rec full}
$\rec(\P)$ is full-dimensional.
\end{claim}

\begin{cpf}
Assume, for a contradiction, that $\rec(\P)$ is not full-dimensional.
\blue{Since $\rec(\P)$ is a rational polyhedron of dimension less than $n$ containing the origin, it is contained in a rational hyperplane through the origin, thus there exists $a \in \Z^n \setminus \bra{\z}$ with $a^\transp d = 0$ for every $d \in \rec(\P)$.
It follows that $a^\transp x$ is bounded on $\P$, since $\P$ is the sum of a polytope and $\rec(\P)$, and it is an integer for every $x \in \P \cap \Z^n$.
Therefore, there exist finitely many integers $c_1,c_2,\dots,c_k$ such that each vector in $\P \cap \Z^n$ lies in one of the rational hyperplanes $\H_i := \bra{x \in \R^n : a^\transp x = c_i}$, for $i \in \bra{1,2,\dots,k}$.}
Claim~\ref{claim H} implies that $f$ is bounded below on $\P \cap \H_i \cap \Z^n$, for every $i \in \bra{1,2,\dots,k}$.
Hence, $f$ is bounded below on $\P \cap \Z^n$, which gives us a contradiction.
\end{cpf}

\begin{claim}
\label{claim origin}
We can assume that the origin $\z$ is in the interior of $\P$.
\end{claim}

\begin{cpf}
It follows from Claim~\ref{claim rec full} that there exists a vector $v \in \Z^n$ in the interior of $\rec(\P)$.
Now let $y \in \P \cap \Z^n$.
It is simple to check that $y+v \in \Z^n$ is in the interior of $\P$.
We can now perform the change of variables $x' = x - (y+v)$, and obtain that $\z$ is in the interior of the resulting polyhedron.
\green{Since $y+v$ is an integer vector, the resulting polyhedron is still rational, pointed, integral, and full-dimensional, it has the same recession cone as $\P$, and it satisfies the property of Claim~\ref{claim H}; moreover, the resulting function is still a polynomial of degree at most three.}
\end{cpf}

\blue{The reductions in Claims~\ref{claim pointed}--\ref{claim origin} modify $\P$.
In the remainder of the proof, $\P = \bra{x \in \R^n : Ax \le b}$ denotes the resulting polyhedron, where $A \in \Z^{m \times n}$ and $b \in \Z^m$, and $a^\transp_1,a^\transp_2,\dots,a^\transp_m$ denote the rows of $A$.}
\green{We remark that $\z \in \P$ and $\rec(\P)$ is a cone imply $\rec(\P) \subseteq \P$, and that $\rec(\P)$ is pointed, since $\P$ is pointed.}

\begin{claim}
\label{claim sequence}
There exists a sequence $\bra{x^j}_{j \in \N}$ of vectors in $\P \cap \Z^n$ satisfying
\begin{subequations}
\label{prop sequence}
\begin{align}
\label{prop lim f}
& \lim_{j \to + \infty} f(x^j) = - \infty, \\
\label{pr f decrease}
& f\pare{x^{j}} > f\pare{x^{j+1}} && \forall j \in \N, \\
\label{prop cone}
& x^j \in \rec(\P) && \forall j \in \N, \\
\label{prop lim norm}
& \lim_{j \to + \infty} \norm{x^j} = \infty, \\
\label{prop norm}
& \norm{x^{j}} < \norm{x^{j+1}} && \forall j \in \N.
\end{align}
\end{subequations}
\end{claim}

\begin{cpf}
Since $f$ is unbounded below on $\P \cap \Z^n$, there exists a sequence $\bra{x^j}_{j \in \N}$ of vectors in $\P \cap \Z^n$ satisfying \eqref{prop lim f}.
By passing to a subsequence of $\bra{x^j}_{j \in \N}$ with monotonically decreasing $f(x^j)$, we can assume that also \eqref{pr f decrease} holds.

Next, we show that there are infinitely many indices $j \in \N$ such that $x^j \in \rec(\P)$.
If not, there exists $k \in \N$ such that, for every $j > k$, $x^j \notin \rec(\P)$.
\blue{For every $j > k$, since $x^j$ is in $\P$ but not in $\rec(\P)$, there exist $i \in \bra{1,2,\dots,m}$ and an integer $c \in \bra{1,\dots,b_i}$ such that $a_i^\transp x^j = c$ (note that $b_i \ge 1$ for every $i$, since $\z$ is in the interior of $\P$).
Since there are finitely many such pairs $(i,c)$, there is at least one hyperplane $\H := \bra{x \in \R^n : a_i^\transp x = c}$ containing infinitely many points among $x^j$, for $j > k$.}
It then follows that $f$ is unbounded below on $\P \cap \H \cap \Z^n$.
This contradicts Claim~\ref{claim H}.

Since there are infinitely many indices $j \in \N$ such that $x^j \in \rec(\P)$, the corresponding subsequence of $\bra{x^j}_{j \in \N}$ also satisfies \eqref{prop cone}.
Property~\eqref{prop lim norm} holds because $f(x^j) \to -\infty$ by~\eqref{prop lim f}, while $f$ is bounded on every bounded set.
By passing to a subsequence of $\bra{x^j}_{j \in \N}$ with monotonically increasing $\norm{x^j}$, we can assume that also \eqref{prop norm} holds.
\end{cpf}

\phead{\blue{Stage 2: restrictions of $f$ to rays.}}
\begin{greenpart}
From Observation~\ref{obs f notation}, we can write $f$, up to an additive constant, which does not affect the statement, in the form
\[
f(x) = T[x,x,x] + M[x,x] + V[x],
\]
where $T \in \R^{n \times n \times n}$ is a symmetric tensor, $M \in \R^{n \times n}$ is a symmetric matrix, and $V \in \R^n$.
For every $y,d \in \R^n$, multilinearity and symmetry give the following identity of polynomials in $\lambda \in \R$:
\begin{align}
\label{eq expansion}
f(y + \lambda d) = T[d,d,d] \, \lambda^3 + \pare{3 T[y,d,d] + M[d,d]} \lambda^2 + \gamma_d(y) \, \lambda + f(y),
\end{align}
where
\begin{align}
\label{eq gamma}
\gamma_d(y) := 3 T[y,y,d] + 2 M[y,d] + V[d].
\end{align}
Recall from Section~\ref{sec notation} that $T[d,v] \in \R^n$ denotes the contraction of $T$ with $d$ and $v$, so that $T[x,d,v] = x^\transp T[d,v]$ for every $x \in \R^n$, and that, by symmetry, $T[x,d,v]$ is invariant under permutations of its arguments.
In particular, the condition $T[d,v] = \z$ means that $T[x,d,v] = 0$ for every $x \in \R^n$.

In the remainder of the proof, we argue by contradiction, and we assume that
\begin{align}
\label{eq star}
\lim_{\nu \to + \infty} f(\bar z + \nu w) \neq - \infty
\qquad \forall \bar z \in \P \cap \Z^n, \ \forall w \in \rec(\P) \setminus \bra{\z}.
\end{align}
Since $\ray(\bar z,w)$ is a ray of $\P$ for every $\bar z \in \P$ and $w \in \rec(\P) \setminus \bra{\z}$, the negation of \eqref{eq star} is precisely the statement of Proposition~\ref{prop int ray}.
We will derive a contradiction from \eqref{eq star}.
Recall that $\z \in \P \cap \Z^n$ by Claim~\ref{claim origin}, thus \eqref{eq star} can be applied with $\bar z = \z$.

\begin{claim}
\label{claim sign}
For every $c \in \R^n \setminus \bra{\z}$, there exists $\bar z \in \P \cap \Z^n$ such that $c^\transp \bar z < 0$.
\end{claim}

\begin{cpf}
Since $\z$ is in the interior of $\P$ (Claim~\ref{claim origin}), there exists $\epsilon > 0$ such that $x^\circ := - \epsilon c \in \P$, and $c^\transp x^\circ = - \epsilon \norm{c}^2 < 0$.
Since $\P$ is integral (Claim~\ref{claim integral}), we have $\P = \conv(\P \cap \Z^n)$, thus $x^\circ = \sum_{l=1}^r \mu_l z_l$ for some $z_1,\dots,z_r \in \P \cap \Z^n$ and $\mu_1,\dots,\mu_r \ge 0$ with $\sum_{l=1}^r \mu_l = 1$.
Since $\sum_{l=1}^r \mu_l \, c^\transp z_l = c^\transp x^\circ < 0$, we obtain $c^\transp z_l < 0$ for some $l \in \bra{1,\dots,r}$.
\end{cpf}

\begin{claim}
\label{claim A}
Let $\bar z \in \P \cap \Z^n$ and $w \in \rec(\P)$.
Then the univariate polynomial $\nu \mapsto f(\bar z + \nu w)$ is either constant or has a positive leading coefficient.
In particular:
\begin{enumerate}[label=(\alph*),ref=(\alph*),noitemsep]
\item
\label{A a}
$T[w,w,w] \ge 0$;
\item
\label{A b}
if $T[w,w,w] = 0$, then $3 T[\bar z,w,w] + M[w,w] \ge 0$;
\item
\label{A c}
if $T[w,w,w] = 0$ and $3 T[\bar z,w,w] + M[w,w] = 0$, then $\gamma_w(\bar z) \ge 0$.
\end{enumerate}
\end{claim}

\begin{cpf}
If $w = \z$, the polynomial is constant.
Otherwise, a nonconstant univariate polynomial with negative leading coefficient tends to $-\infty$ as $\nu \to + \infty$, which is excluded by \eqref{eq star}.
Items~\ref{A a}--\ref{A c} follow by reading off the coefficients in \eqref{eq expansion} with $y = \bar z$ and $d = w$.
\end{cpf}

For every unit vector $e \in \rec(\P)$, we define
\[
\ell(e) := \inf \bra{f(\lambda e) : \lambda \ge 0}.
\]
By Claim~\ref{claim A} with $\bar z = \z$, the polynomial $\lambda \mapsto f(\lambda e)$ is bounded below on $\R_{\ge 0}$, thus $\ell(e) \in \R$.
The next two claims are the analytic core of the proof.
In both claims, $\bra{e^j}_{j \in \N}$ denotes a sequence of unit vectors in $\rec(\P)$ such that $\ell(e^j) \to -\infty$, as $j \to +\infty$.
We use the shorthand $T_j := T[e^j,e^j,e^j]$, $M_j := M[e^j,e^j]$, $V_j := V[e^j]$, so that
\begin{align}
\label{eq gj}
f(\lambda e^j) = T_j \lambda^3 + M_j \lambda^2 + V_j \lambda,
\end{align}
and we note that $T_j \ge 0$ by Claim~\ref{claim A}\ref{A a}, and that $\abs{M_j} \le \norm{M}$ and $\abs{V_j} \le \norm{V}$ for every $j \in \N$.

\begin{claim}
\label{claim B}
Let $\bra{e^j}_{j \in \N}$ be a sequence of unit vectors in $\rec(\P)$ with $\ell(e^j) \to -\infty$, and assume that $e^j \to e$.
Then $T[e,e,e] = 0$, $M[e,e] = 0$, and $V[e] \ge 0$.
\end{claim}

\begin{cpf}
Note that $e \in \rec(\P)$, since $\rec(\P)$ is closed, and that $T_j \to T[e,e,e]$, $M_j \to M[e,e]$, $V_j \to V[e]$ by continuity.

Since $T_j \ge 0$ for every $j$, we have $T[e,e,e] \ge 0$.
Assume, for a contradiction, $T[e,e,e] > 0$, and let $\delta := T[e,e,e]/2$.
For $j$ large enough we have $T_j \ge \delta$, thus, for every $\lambda \ge 0$,
\[
f(\lambda e^j) \ge \delta \lambda^3 - \norm{M} \lambda^2 - \norm{V} \lambda.
\]
The right-hand side is a cubic with positive leading coefficient, thus it is bounded below on $\R_{\ge 0}$ by a constant independent of $j$.
Hence $\ell(e^j)$ is bounded below, contradicting $\ell(e^j) \to -\infty$.
This shows $T[e,e,e] = 0$.

Since $T[e,e,e] = 0$, Claim~\ref{claim A}\ref{A b} with $w = e$ and $\bar z = \z$ gives $M[e,e] \ge 0$.
Assume, for a contradiction, $M[e,e] > 0$.
For $j$ large enough we have $M_j \ge M[e,e]/2$, thus, for every $\lambda \ge 0$, using $T_j \lambda^3 \ge 0$,
\[
f(\lambda e^j) \ge \frac{M[e,e]}{2} \lambda^2 - \norm{V} \lambda \ge - \frac{\norm{V}^2}{2 M[e,e]}.
\]
Hence $\ell(e^j)$ is bounded below, contradicting $\ell(e^j) \to -\infty$.
This shows $M[e,e] = 0$.

Finally, Claim~\ref{claim A}\ref{A c} with $w = e$ and $\bar z = \z$ gives $V[e] = \gamma_e(\z) \ge 0$.
\end{cpf}

\begin{claim}
\label{claim C}
Let $\bra{e^j}_{j \in \N}$ be a sequence of unit vectors in $\rec(\P)$ with $\ell(e^j) \to -\infty$.
Then it is impossible that both of the following conditions hold:
\begin{enumerate}[label=(\roman*),ref=(\roman*),noitemsep]
\item
\label{C i}
$M_j > 0$ for every $j$ large enough;
\item
\label{C ii}
there exists a constant $K \ge 0$ such that $\max\bra{0, -V_j} \le K M_j$ for every $j$ large enough.
\end{enumerate}
\end{claim}

\begin{cpf}
Assume that both~\ref{C i} and~\ref{C ii} hold.
Let $g_j(\lambda) := f(\lambda e^j)$, as in \eqref{eq gj}.
For every $j$ large enough and every $\lambda \ge K$, we have
\[
g_j'(\lambda) = 3 T_j \lambda^2 + 2 M_j \lambda + V_j \ge 2 M_j \lambda + V_j \ge 2 K M_j + V_j \ge 2 \max\bra{0,-V_j} + V_j \ge 0,
\]
where we used $T_j \ge 0$, $M_j > 0$, and~\ref{C ii}.
Hence $g_j$ is nondecreasing on $[K,+\infty)$, and
\[
\ell(e^j) = \min \bra{g_j(\lambda) : \lambda \in [0,K]} \ge \min \bra{f(x) : \norm{x} \le K},
\]
which is a finite constant independent of $j$, since $f$ is continuous.
This contradicts $\ell(e^j) \to -\infty$.
\end{cpf}

\phead{\blue{Stage 3: }chains of projections.}
We now introduce the objects that are iteratively constructed in the third stage of the proof.
Given a face $\F$ of $\rec(\P)$ and a vector $u \in \F \setminus \bra{\z}$, we will project points of $\F$ along $-u$ onto a proper face of $\F$.
The properties of the resulting sequences are collected in the following definition.
Let $k \ge 0$ be an integer.
A \emph{chain of order $k$} consists of
\begin{itemize}[noitemsep]
\item
faces $\rec(\P) = \F_0 \supsetneq \F_1 \supsetneq \dots \supsetneq \F_k$ of $\rec(\P)$;
\item
unit vectors $d_1,\dots,d_k$ with $d_i \in \F_{i-1}$, for $i \in \bra{1,\dots,k}$;
\item
for each $i \in \bra{0,1,\dots,k}$, a sequence $\bra{x_i^j}_{j \in \N}$ of vectors in $\F_i \setminus \bra{\z}$ with $f(x_i^j) \to -\infty$, as $j \to +\infty$, such that $x_0^j \in \Z^n$ for every $j \in \N$, and, for $i \ge 1$, $x_i^j = x_{i-1}^j - p_{i,j} d_i$ for some $p_{i,j} \ge 0$;
\end{itemize}
such that, writing
\[
e_i^j := \frac{x_{i-1}^j}{\norm{x_{i-1}^j}} \qquad \forall i \in \bra{1,\dots,k+1}, \ \forall j \in \N,
\]
the following conditions hold:
\begin{enumerate}[label=(C\arabic*),ref=(C\arabic*),noitemsep,leftmargin=3em]
\item
\label{C1}
$e_i^j \to d_i$, as $j \to +\infty$, for every $i \in \bra{1,\dots,k}$;
\item
\label{C2}
for every $i \in \bra{1,\dots,k}$, there exist $\sigma_{i,j} \in \R$ and $\tau_{i,j} > 0$, for $j \in \N$, such that $e_i^j = \sigma_{i,j} d_i + \tau_{i,j} e_{i+1}^j$ for every $j \in \N$, and $\sigma_{i,j} \to 1$, $\tau_{i,j} \to 0$, as $j \to +\infty$;
\item
\label{C3}
$T[d_a,d_b] = \z$ for every $a,b \in \bra{1,\dots,k}$;
\item
\label{C4}
$M[d_a,d_b] = 0$ for every $a,b \in \bra{1,\dots,k}$;
\item
\label{C5}
$V[d_a] \ge 0$ for every $a \in \bra{1,\dots,k}$;
\item
\label{C6}
$\gamma_{d_a}(\bar z) \ge 0$ for every $a \in \bra{1,\dots,k}$ and every $\bar z \in \P \cap \Z^n$.
\end{enumerate}
We make some remarks.
The vectors $x_i^j$ with $i \ge 1$ are, in general, not integer; only the vectors $x_0^j$ are required to be integer.
Since $\F_i \subseteq \rec(\P) \subseteq \P$, all vectors $x_i^j$ lie in $\P$, and $\ray(\z,e_i^j)$ is a ray of $\P$.
Since $f$ is bounded on every bounded set, $f(x_i^j) \to -\infty$ implies $\norm{x_i^j} \to \infty$.
All conditions in the definition are preserved if we replace $j$ by a subsequence, simultaneously for all $i \in \bra{0,\dots,k}$; we will do so several times without further comment.
Since the inclusions $\F_{i-1} \supsetneq \F_i$ are strict, we have $\dim(\F_k) \le \dim(\rec(\P)) - k$.
For $k=0$, conditions~\ref{C1}--\ref{C6} are vacuous, and the sequence $\bra{x^j}_{j \in \N}$ of Claim~\ref{claim sequence} (after removing the vector $\z$, which can appear at most once) is a chain of order~$0$, with $x_0^j := x^j$ and $\F_0 := \rec(\P)$.

\begin{claim}
\label{claim E}
Let $d_1,\dots,d_k \in \R^n$ satisfy \ref{C3}, \ref{C4}, and \ref{C6}.
Let $\bar z \in \P \cap \Z^n$ and $w := \sum_{a=1}^k \lambda_a d_a$ with $\lambda_1,\dots,\lambda_k \ge 0$.
Then, for every $\nu \in \R$,
\[
f(\bar z + \nu w) = f(\bar z) + \nu \sum_{a=1}^k \lambda_a \gamma_{d_a}(\bar z),
\qquad \text{and} \qquad
\sum_{a=1}^k \lambda_a \gamma_{d_a}(\bar z) \ge 0.
\]
In particular, $f(\bar z - w) \le f(\bar z)$.
\end{claim}

\begin{cpf}
We use \eqref{eq expansion} with $y = \bar z$ and $d = w$.
By multilinearity and symmetry,
\begin{align*}
T[w,w,w] & = \sum_{a,b,c=1}^k \lambda_a \lambda_b \lambda_c \, d_c^\transp T[d_a,d_b] = 0, \\
3 T[\bar z,w,w] + M[w,w] & = \sum_{a,b=1}^k \lambda_a \lambda_b \pare{3 \, \bar z^\transp T[d_a,d_b] + M[d_a,d_b]} = 0,
\end{align*}
by~\ref{C3} and~\ref{C4}.
Furthermore, $\gamma_w(\bar z)$ is linear in $w$, thus $\gamma_w(\bar z) = \sum_{a=1}^k \lambda_a \gamma_{d_a}(\bar z) \ge 0$ by~\ref{C6}.
Since \eqref{eq expansion} is an identity of polynomials in $\nu$, it holds for every real $\nu$, in particular for $\nu = -1$.
\end{cpf}

\begin{claim}
\label{claim dim}
If a chain of order $k$ exists, then $\dim(\F_k) \ge 2$.
In particular, $k \le \dim(\rec(\P)) - 2$.
\end{claim}

\begin{cpf}
Since $\F_k$ contains the nonzero vectors $x_k^j$, we have $\dim(\F_k) \ge 1$.
Assume $\dim(\F_k) = 1$.
Since $\F_k$ is a face of the pointed rational cone $\rec(\P)$, it is a pointed rational cone of dimension one, thus $\F_k = \ray(\z,u)$ for some $u \in \Z^n \setminus \bra{\z}$.
Hence $x_k^j = \mu_j u$ with $\mu_j > 0$, and $\mu_j \to +\infty$ because $\norm{x_k^j} \to \infty$.
Since $f(\mu_j u) \to -\infty$ along a sequence $\mu_j \to +\infty$, the univariate polynomial $\mu \mapsto f(\mu u)$ tends to $-\infty$ as $\mu \to +\infty$.
This contradicts \eqref{eq star} with $\bar z = \z$ and $w = u$.
The last statement follows from $\dim(\F_k) \le \dim(\rec(\P)) - k$.
\end{cpf}

The next claim is the heart of the proof.

\begin{claim}
\label{claim F}
If a chain of order $k$ exists, then a chain of order $k+1$ exists.
\end{claim}

\begin{cpf}
Fix a chain of order $k$, with the notation introduced above.
Let $e^j := e_{k+1}^j = x_k^j / \norm{x_k^j}$.
The vectors $e^j$ are unit vectors in the closed cone $\F_k$, thus, by passing to a subsequence, we can assume $e^j \to u$, for some unit vector $u \in \F_k \subseteq \rec(\P)$.
The proof is subdivided into ten steps.

\phead{Step 1: $T[u,u,u] = 0$, $M[u,u] = 0$, and $V[u] \ge 0$.}
Since $x_k^j = \norm{x_k^j} e^j$ lies on $\ray(\z,e^j)$, we have $\ell(e^j) \le f(x_k^j) \to -\infty$.
The statement then follows from Claim~\ref{claim B}.

\phead{Step 2: $T[u,u] = \z$.}
Assume, for a contradiction, $T[u,u] \neq \z$.
By Claim~\ref{claim sign}, there exists $\bar z \in \P \cap \Z^n$ with $T[\bar z,u,u] = \bar z^\transp T[u,u] < 0$.
In \eqref{eq expansion} with $y = \bar z$ and $d = u$, the cubic coefficient is $T[u,u,u] = 0$ by Step~1, and the quadratic coefficient is $3 T[\bar z,u,u] + M[u,u] = 3 T[\bar z,u,u] < 0$.
This contradicts Claim~\ref{claim A}\ref{A b}.

\phead{Step 3: $\gamma_u(\bar z) \ge 0$ for every $\bar z \in \P \cap \Z^n$.}
By Steps~1 and~2, for every $\bar z \in \R^n$, the cubic and quadratic coefficients in \eqref{eq expansion} with $y = \bar z$ and $d = u$ vanish, since $3 T[\bar z,u,u] + M[u,u] = 3 \, \bar z^\transp T[u,u] = 0$.
The statement then follows from Claim~\ref{claim A}\ref{A c}.

\phead{Step 4: $M[d_a,u] \ge 0$ for every $a \in \bra{1,\dots,k}$.}
Let $w := d_a + u$, which lies in $\rec(\P)$.
We have $w \neq \z$, since otherwise both $u$ and $-u = d_a$ would lie in $\rec(\P)$, contradicting the fact that $\rec(\P)$ is pointed.
By multilinearity and symmetry,
\begin{align}
\label{eq Tw}
T[w,w,w] = d_a^\transp T[d_a,d_a] + 3 \, u^\transp T[d_a,d_a] + 3 \, d_a^\transp T[u,u] + u^\transp T[u,u] = 0,
\end{align}
by~\ref{C3} and Step~2.
Furthermore, $M[w,w] = M[d_a,d_a] + 2 M[d_a,u] + M[u,u] = 2 M[d_a,u]$ by~\ref{C4} and Step~1.
Claim~\ref{claim A}\ref{A b} with $\bar z = \z$ gives $M[d_a,u] \ge 0$.

\phead{Step 5: $M[d_a,u] = 0$ for every $a \in \bra{1,\dots,k}$.}
Fix $a \in \bra{1,\dots,k}$ and assume, for a contradiction, $M[d_a,u] > 0$.
We apply Claim~\ref{claim C} to the sequence $\bra{e_a^j}_{j \in \N}$ of unit vectors in $\F_{a-1} \subseteq \rec(\P)$.
Since $x_{a-1}^j$ lies on $\ray(\z,e_a^j)$, we have $\ell(e_a^j) \le f(x_{a-1}^j) \to -\infty$.

By applying~\ref{C2} for $i = a, a+1, \dots, k$, we obtain, for every $j \in \N$,
\begin{align}
\label{eq unroll}
e_a^j = \sum_{b=a}^k c_{b,j} \, d_b + \theta_j \, e^j,
\end{align}
where
\[
c_{a,j} := \sigma_{a,j}, \qquad
c_{b,j} := \sigma_{b,j} \prod_{i=a}^{b-1} \tau_{i,j} \quad \forall b \in \bra{a+1,\dots,k}, \qquad
\theta_j := \prod_{i=a}^{k} \tau_{i,j}.
\]
Indeed, $e_a^j = \sigma_{a,j} d_a + \tau_{a,j} e_{a+1}^j$, $e_{a+1}^j = \sigma_{a+1,j} d_{a+1} + \tau_{a+1,j} e_{a+2}^j$, and so on until $e_k^j = \sigma_{k,j} d_k + \tau_{k,j} e_{k+1}^j = \sigma_{k,j} d_k + \tau_{k,j} e^j$.
From~\ref{C2} we obtain: $\theta_j > 0$ for every $j$; $c_{a,j} \to 1$; $c_{b,j} \to 0$ for every $b \in \bra{a+1,\dots,k}$; $\theta_j \to 0$; $c_{b,j} > 0$ for every $b \in \bra{a,\dots,k}$ and $j$ large enough; and, for every $b \in \bra{a,\dots,k}$,
\begin{align}
\label{eq ratio}
\frac{\theta_j}{c_{b,j}} = \frac{\prod_{i=b}^k \tau_{i,j}}{\sigma_{b,j}} \longrightarrow 0.
\end{align}
We now use \eqref{eq unroll} to expand the coefficients $M_j = M[e_a^j,e_a^j]$ and $V_j = V[e_a^j]$ in \eqref{eq gj}, where $\bra{e_a^j}_{j \in \N}$ plays the role of $\bra{e^j}_{j \in \N}$.

By~\ref{C4}, all the terms $M[d_b,d_{b'}]$ vanish, thus
\[
M[e_a^j,e_a^j] = 2 \theta_j \sum_{b=a}^k c_{b,j} M[d_b,e^j] + \theta_j^2 M[e^j,e^j] = \theta_j \, m_j,
\]
where
\[
m_j := 2 \sum_{b=a}^k c_{b,j} M[d_b,e^j] + \theta_j M[e^j,e^j].
\]
Since $c_{a,j} \to 1$, $c_{b,j} \to 0$ for $b > a$, $\theta_j \to 0$, $e^j \to u$, and $\abs{M[x,y]} \le \norm{M}$ for all unit vectors $x,y$, we obtain $m_j \to 2 M[d_a,u] > 0$.
Hence $M[e_a^j,e_a^j] = \theta_j m_j > 0$ for $j$ large enough, which is condition~\ref{C i} of Claim~\ref{claim C}.

By linearity, $V[e_a^j] = \sum_{b=a}^k c_{b,j} V[d_b] + \theta_j V[e^j]$.
We distinguish two cases.

\phead{Case (a): $V[d_{b^*}] > 0$ for some $b^* \in \bra{a,\dots,k}$.}
For $j$ large enough all $c_{b,j}$ are positive, and $V[d_b] \ge 0$ by~\ref{C5}, thus
\[
V[e_a^j] \ge c_{b^*,j} V[d_{b^*}] - \theta_j \norm{V} = c_{b^*,j} \pare{ V[d_{b^*}] - \frac{\theta_j}{c_{b^*,j}} \norm{V} } > 0
\]
for $j$ large enough, by \eqref{eq ratio}.
Hence $\max\bra{0,-V[e_a^j]} = 0$ for $j$ large enough, and condition~\ref{C ii} of Claim~\ref{claim C} holds with $K = 0$.

\phead{Case (b): $V[d_b] = 0$ for every $b \in \bra{a,\dots,k}$.}
Then $V[e_a^j] = \theta_j V[e^j]$, thus $\abs{V[e_a^j]} \le \theta_j \norm{V}$.
Since $m_j \to 2 M[d_a,u]$, for $j$ large enough we have $m_j \ge M[d_a,u]$, thus
\[
\max\bra{0,-V[e_a^j]} \le \theta_j \norm{V} = \frac{\norm{V}}{m_j} \, M[e_a^j,e_a^j] \le \frac{\norm{V}}{M[d_a,u]} \, M[e_a^j,e_a^j],
\]
and condition~\ref{C ii} of Claim~\ref{claim C} holds with $K = \norm{V}/M[d_a,u]$.

In both cases, conditions~\ref{C i} and~\ref{C ii} hold, contradicting Claim~\ref{claim C}.
Hence $M[d_a,u] = 0$.

\phead{Step 6: $T[d_a,u] = \z$ for every $a \in \bra{1,\dots,k}$.}
Assume, for a contradiction, $T[d_a,u] \neq \z$ for some $a$.
By Claim~\ref{claim sign}, there exists $\bar z \in \P \cap \Z^n$ with $T[\bar z,d_a,u] = \bar z^\transp T[d_a,u] < 0$.
Let $w := d_a + u \in \rec(\P) \setminus \bra{\z}$, as in Step~4.
By \eqref{eq Tw}, the cubic coefficient in \eqref{eq expansion} with $y = \bar z$ and $d = w$ vanishes, while the quadratic coefficient is
\begin{align*}
3 T[\bar z,w,w] + M[w,w]
& = 3 \pare{\bar z^\transp T[d_a,d_a] + 2 \, \bar z^\transp T[d_a,u] + \bar z^\transp T[u,u]} \\
& \quad + \pare{M[d_a,d_a] + 2 M[d_a,u] + M[u,u]} \\
& = 6 \, T[\bar z,d_a,u] < 0,
\end{align*}
by~\ref{C3}, Step~2, \ref{C4}, Step~5, and Step~1.
This contradicts Claim~\ref{claim A}\ref{A b}.

\phead{Step 7: the new direction.}
Let $d_{k+1} := u \in \F_k$.
Steps~1, 2, 5, and~6, together with~\ref{C3} and~\ref{C4}, give $T[d_a,d_b] = \z$ and $M[d_a,d_b] = 0$ for every $a,b \in \bra{1,\dots,k+1}$.
Step~1 and~\ref{C5} give $V[d_a] \ge 0$, and Step~3 and~\ref{C6} give $\gamma_{d_a}(\bar z) \ge 0$ for every $\bar z \in \P \cap \Z^n$, for every $a \in \bra{1,\dots,k+1}$.
That is, conditions~\ref{C3}--\ref{C6} hold for $d_1,\dots,d_{k+1}$.
In particular, Claim~\ref{claim E} applies to $d_1,\dots,d_{k+1}$.
Furthermore, $e_{k+1}^j = e^j \to u = d_{k+1}$, which is condition~\ref{C1} for $i = k+1$.

\phead{Step 8: the projection.}
Since $u \in \rec(\P) \setminus \bra{\z}$ and $\rec(\P)$ is pointed, we have $-u \notin \rec(\P)$, thus the index set $I^- := \bra{i \in \bra{1,\dots,m} : a_i^\transp u < 0}$ is nonempty.
For every $j \in \N$, we define
\[
q_j := \min \bra{\frac{a_i^\transp x_k^j}{a_i^\transp u} : i \in I^-} \ge 0,
\qquad
x_{k+1}^j := x_k^j - q_j u,
\]
where $q_j \ge 0$ because $a_i^\transp x_k^j \le 0$ and $a_i^\transp u < 0$ for every $i \in I^-$.
We show that $x_{k+1}^j \in \F_k$.
Since $\F_k$ is a face of $\rec(\P)$, we can write $\F_k = \rec(\P) \cap \bra{x \in \R^n : a_i^\transp x = 0 \ \forall i \in I_k}$ for some $I_k \subseteq \bra{1,\dots,m}$.
For $i \in I^-$, the definition of $q_j$ gives $a_i^\transp x_{k+1}^j = a_i^\transp x_k^j - q_j a_i^\transp u \le 0$.
For $i \notin I^-$, we have $a_i^\transp u \ge 0$, thus $a_i^\transp x_{k+1}^j \le a_i^\transp x_k^j \le 0$.
Hence $x_{k+1}^j \in \rec(\P)$.
For $i \in I_k$, we have $a_i^\transp u = 0$ because $u \in \F_k$, thus $a_i^\transp x_{k+1}^j = a_i^\transp x_k^j = 0$.
Hence $x_{k+1}^j \in \F_k$.

Let $i_j \in I^-$ be an index attaining the minimum in the definition of $q_j$, so that $a_{i_j}^\transp x_{k+1}^j = 0$.
Since $I^-$ is finite, by passing to a subsequence we can assume that $i_j = i^*$ for every $j \in \N$.
We define
\[
\F_{k+1} := \F_k \cap \bra{x \in \R^n : a_{i^*}^\transp x = 0}.
\]
Since the inequality $a_{i^*}^\transp x \le 0$ is valid for $\rec(\P) \supseteq \F_k$, the set $\F_{k+1}$ is a face of $\F_k$, and therefore a face of $\rec(\P)$.
It is a proper face of $\F_k$, because $u \in \F_k$ and $a_{i^*}^\transp u < 0$.
By construction, $x_{k+1}^j \in \F_{k+1}$ for every $j \in \N$.

\phead{Step 9: $f(x_{k+1}^j) \to -\infty$.}
By the definition of a chain, $x_0^j = x_{k+1}^j + w^j$, where
\[
w^j := \sum_{a=1}^k p_{a,j} d_a + q_j d_{k+1}, \qquad p_{1,j},\dots,p_{k,j},q_j \ge 0.
\]
Since $x_0^j \in \P \cap \Z^n$, Claim~\ref{claim E} applied to $d_1,\dots,d_{k+1}$ (see Step~7) yields $f(x_{k+1}^j) = f(x_0^j - w^j) \le f(x_0^j) \to -\infty$.
Since $f$ is bounded on every bounded set, $\norm{x_{k+1}^j} \to \infty$, thus, by discarding finitely many indices, we can assume $x_{k+1}^j \neq \z$ for every $j \in \N$.

\phead{Step 10: condition~\ref{C2} for $i = k+1$.}
Let $e_{k+2}^j := x_{k+1}^j / \norm{x_{k+1}^j}$, and
\[
\sigma_{k+1,j} := \frac{q_j}{\norm{x_k^j}}, \qquad
\tau_{k+1,j} := \frac{\norm{x_{k+1}^j}}{\norm{x_k^j}} > 0.
\]
Then
\[
e_{k+1}^j = \frac{x_k^j}{\norm{x_k^j}} = \frac{x_{k+1}^j + q_j u}{\norm{x_k^j}} = \sigma_{k+1,j} \, d_{k+1} + \tau_{k+1,j} \, e_{k+2}^j.
\]
From $a_{i^*}^\transp x_{k+1}^j = 0$ we obtain $q_j = a_{i^*}^\transp x_k^j / a_{i^*}^\transp u$, hence
\[
\sigma_{k+1,j} = \frac{a_{i^*}^\transp e_{k+1}^j}{a_{i^*}^\transp u} \longrightarrow \frac{a_{i^*}^\transp u}{a_{i^*}^\transp u} = 1,
\qquad
\tau_{k+1,j} = \norm{e_{k+1}^j - \sigma_{k+1,j} u} \longrightarrow \norm{u - u} = 0.
\]

\phead{Conclusion.}
The faces $\F_0 \supsetneq \dots \supsetneq \F_k \supsetneq \F_{k+1}$, the unit vectors $d_1,\dots,d_{k+1}$ with $d_{k+1} = u \in \F_k$, and the sequences $\bra{x_i^j}_{j \in \N}$, for $i \in \bra{0,\dots,k+1}$, with $x_{k+1}^j = x_k^j - q_j d_{k+1}$, $x_{k+1}^j \in \F_{k+1} \setminus \bra{\z}$, and $f(x_{k+1}^j) \to -\infty$ (Steps~8 and~9), satisfy conditions~\ref{C1} and~\ref{C2} for every $i \in \bra{1,\dots,k+1}$ (by the given chain, Step~7, and Step~10), and conditions~\ref{C3}--\ref{C6} for every $a,b \in \bra{1,\dots,k+1}$ (Step~7).
This is a chain of order $k+1$.
\end{cpf}
\end{greenpart}

\green{We can now conclude the proof of Proposition~\ref{prop int ray}.
As observed above, Claim~\ref{claim sequence} provides a chain of order~$0$.
By Claim~\ref{claim F} and induction on $k$, a chain of order $k$ exists for every $k \ge 0$.
In particular, a chain of order $n-1$ exists, and it satisfies $\dim(\F_{n-1}) \le \dim(\rec(\P)) - (n-1) \le n - (n-1) = 1$, contradicting Claim~\ref{claim dim}.
Hence \eqref{eq star} fails: there exist $\bar z \in \P \cap \Z^n$ and $w \in \rec(\P) \setminus \bra{\z}$ such that $f$ is unbounded below on the ray $\ray(\bar z,w)$ of $\P$.}
\end{prfc}

\ifthenelse {\boolean{MOR}}
{
}
{

\bigskip
\noindent
\textbf{Acknowledgments:} The author thanks Daniel Bienstock and Robert Hildebrand for valuable discussions on unbounded rays.

\bigskip
\noindent
\textbf{Funding:}
Alberto~Del~Pia is partially funded by AFOSR grant FA9550-23-1-0433 and ONR grant N00014-25-1-2490. Any opinions, findings, and conclusions or recommendations expressed in this material are those of the author and do not necessarily reflect the views of the Air Force Office of Scientific Research or of the Office of Naval Research.

}

\end{document}